 \documentclass[psamsfonts]{article}
 
\usepackage{amsfonts}
\usepackage{amsmath}
\usepackage{amssymb}
\usepackage{amscd}
\usepackage{amstext}
\usepackage{amsthm}

 \title{Extensions by simple $C^*$-algebras\\
       -- Quasidiagonal extensions
 \thanks{Research partially supported by NSF grants DMS 0097903.
         AMS 1991 Subject Classification Numbers:
                         Primary 46L05,
                                46L35.
                        Key words: Extensions, Simple $C^*$-algebras.
                                                        \protect\\}}
\author{{\bf Dedicated to Lawrence G. Brown on his 60th birthday}\\
Huaxin Lin\\
Department of Mathematics\\
University of Oregon\\
Eugene, Oregon 97403-1222}
\date{}
\begin{document}
\maketitle

\newcommand{\CA}{$C^*$-algebra}
\newcommand{\SCA}{$C^*$-subalgebra}
\newcommand{\aue}{approximate unitary equivalence}
\newcommand{\ayue}{approximately unitarily equivalent}
\newcommand{\mops}{mutually orthogonal projections}
\newcommand{\hm}{homomorphism}
\newcommand{\pisca}{purely infinite simple \CA}
\newcommand{\andeqn}{\,\,\,\,\,\, {\rm and} \,\,\,\,\,\,}
\newcommand{\QED}{\rule{1.5mm}{3mm}}
\newcommand{\morp}{contractive completely
positive linear
 map}
\newcommand{\asmorp}{asymptotic morphism}
\newcommand{\arrow}{\rightarrow}
\newcommand{\tdsum}{\widetilde{\oplus}}
\newcommand{\pa}{\|}  
\newcommand{\ep}{\varepsilon}
\newcommand{\id}{{\rm id}}
\newcommand{\aueeps}[1]{\stackrel{#1}{\sim}}
\newcommand{\aeps}[1]{\stackrel{#1}{\approx}}
\newcommand{\dt}{\delta}
\newcommand{\yu}{\fang}
\newcommand{\ca}{{\cal C}_1}
\newcommand{\Ad}{{\rm ad}}

\newtheorem{thm}{Theorem}[section]
\newtheorem{Lem}[thm]{Lemma}
\newtheorem{Prop}[thm]{Proposition}
\newtheorem{Def}[thm]{Definition}
\newtheorem{Cor}[thm]{Corollary}
\newtheorem{Ex}[thm]{Example}
\newtheorem{Pro}[thm]{Problem}
\newtheorem{Remark}[thm]{Remark}
\newtheorem{NN}[thm]{}
\renewcommand{\theequation}{e\,\arabic{section}.\arabic{equation}}

\newcommand{\Ik}{ {\cal I}^{(k)}}
\newcommand{\Iz}{{\cal I}^{(0)}}
\newcommand{\Ii}{{\cal I}^{(1)}}
\newcommand{\Ip}{{\cal I}^{(2)}}

\begin{abstract}
Let $A$ be an amenable separable \CA\,  and $B$ be a non-unital
but $\sigma$-unital simple \CA\, with continuous scale.
We show that two essential extensions
$\tau_1$ and $\tau_2$ of $A$ by $B$ are  approximately
unitarily equivalent if and only if
$$
[\tau_1]=[\tau_2]\,\,\, {\rm in}\,\,\, KL(A, M(B)/B).
$$
If $A$ is assumed to satisfy the Universal Coefficient Theorem,
there is a bijection from approximate unitary equivalence
classes of the above mentioned extensions to
$KL(A, M(B)/B).$
Using $KL(A, M(B)/B),$ we compute exactly when an essential extension
is quasidiagonal. We show that quasidiagonal extensions
may not be approximately trivial.
We also study the approximately trivial extensions.
\end{abstract}

\vspace{0.2in}
\noindent

{\bf \Large{ Introduction}}
\vspace{0.1in}

The study of $C^*$-algebra extensions of $C(X)$ by compact operators
was motivated by the understanding of essentially normal operators
on an infinite dimensional Hilbert space.
The Brown-Douglas-Fillmore Theory for essentially
normal operators gives the classification of essentially normal operators
up to unitary equivalence (\cite{BDF1}). The original
BDF-theory quickly developed
into   \CA\, extension theory (\cite{BDF2} and \cite{BDF3})
and the KK-theory of Kasparov.
Applications of this development can be found not only in
operator theory and operator algebras but also in both
geometry and non-commutative geometry.

Let $0\to B\to E\to A\to 0$ be an essential extension of $A$ by $B.$
This is determined by a monomorphism $\tau: A\to M(B)/B.$
While $KK$-theory gives the classification of extensions up to stable
unitary equivalence, it does not give much information
on essential extensions when $B\not ={\cal K},$ where ${\cal K}$ is the
compact operators on $l^2.$
The example in 1.1 below shows that a non-trivial extension $\tau$
may have $[\tau]=0$ in $KK^1(A, B).$ Other examples (see \ref{IIEx}) show
that there may be infinitely many non-equivalent trivial extensions.
Extensions by simple \CA s have been studied in a few special cases
(see \cite{Ln3}, \cite{Ln4},
\cite{Ln5} and \cite{Lnamj}).

  In this paper, we study approximately unitary equivalence
classes of essential extensions of separable amenable \CA s by
$\sigma$-unital simple \CA s.
One of the reasons that BDF-theory was successful  is that
the Calkin algebra $M({\cal K})/{\cal K}$ is simple (and purely infinite).
We will restrict ourselves to the case that $M(B)/B$ is simple.
It is shown in \cite{Ln1} and more recently in \cite{Lnnew} that,
for a non-unital and $\sigma$-unital simple \CA\, $A\not\cong {\cal K},$
$M(A)/A$ is simple if and only if $A$ has a continuous scale.
Furthermore, in \cite{Lnnew} it is shown that when $M(A)/A$ is simple
it is necessarily  purely infinite simple.

With the Busby invariant, to study essential extensions of $A$ by $B$
it is sufficient to  study monomorphisms from $A$ to $M(B)/B.$
With recent development in classification of
simple amenable \CA s, we know a great deals concerning monomorphisms from
one amenable (simple) \CA\, to a separable  amenable  purely infinite
simple \CA\, (see for example, \cite{Ph1}, \cite{Lnpro} and \cite{Lnsemi}).
However, $M(B)/B$ is not amenable and we do not assume
that $A$ is simple.
One of the main results of this article is the following:
Two essential extensions are approximately
unitarily equivalent if they induce the same element in $KL(A, M(B)/B).$
If furthermore $A$ satisfies the Universal Coefficient
Theorem, then
there is a bijection between
${\bf Ext}_{ap}(A,B),$ the  approximate unitary equivalence
classes of essential extensions, and $KL(A, M(B)/B).$

However, unlike the classical case, the zero element in $KL(A,M(B)/B)$
does not in general give an approximately trivial extension.
On the contrary, at least in some cases, $[\tau]=0$ in $KL(A,M(B)/B)$
never gives an approximately trivial extension and only when $[\tau]\not=0$
in $KL(A, M(B)/B)$ may approximately trivial extensions  occur.
To make matters worse, there may not be any essential trivial extensions
of $A$ by $B$ even though we can use the above mentioned bijection
to classify extensions. This leads us to study quasidiagonal extensions.

Quasidiagonality was  defined by P. R. Halmos (\cite{H}) in 1970.
The \CA\, version soon appeared. L. G. Brown, R. G. Douglas and P. A.
Fillmore (\cite{BDF2}) first recognized that the study of quasidiagonal
extensions might be approached by $K$-theory.
They noticed  that limits of trivial extensions correspond to
the quasidiagonal extensions. L. G. Brown pursued this further  in
(see \cite{Br2}). Further developments in the study of quasidiagonality
can be found in \cite{Sa}, \cite{V1} and \cite{V2}. C. L. Schochet
proved that stable quasidiagonal extensions are the same as
limits of stable trivial extensions and can be characterized by
$Pext(K_*(A), K_*(B))$ if $A$ is assumed to be quasidiagonal relative to $B$
and it satisfies the Universal Coefficient Theorem.
These results might lead one to believe that quasidiagonal extensions
are the same as limits of trivial extensions in  greater generality.
However, in this paper we show this fails when $B$ is neither
isomorphic to ${\cal K}$ nor
purely infinite simple.

One should note that the existence of  quasidiagonal extensions
implies that $B$ has at least one approximate identity consisting
of projections. Suppose that $B$ is a non-unital and
$\sigma$-unital simple \CA\, with real rank zero, stable rank one
and weakly unperforated $K_0(B).$ If $A$ is a quasidiagonal \CA,
then there exists an essential quasidiagonal extension of $A$ by
$B.$ This condition is  necessary if we  assume that $B$ is also a
quasidiagonal \CA. Using $K$-theory and the classification result
mentioned above, we give a necessary and sufficient condition for
essential extensions to be quasidiagonal for a large class of
amenable quasidiagonal \CA s $A.$ We also give a necessary
condition for essential extensions to be approximately trivial for
amenable \CA s which satisfy the UCT. As a consequence, a large
class  of quasidiagonal extensions are {\it not} the limits of
trivial extensions.

The essential extensions of a separable amenable \CA\, $A$ by $B$
(where $B$ is a non-unital and $\sigma$-unital \CA \, with a continuous scale)
is proved in this paper to be determined by $KL(A, M(B)/B).$ However,
to determine which elements in $KL(A, M(B)/B)$ give an approximately trivial
extension is still a difficult task. As mentioned above, for example,
the zero element in $KL(A,M(B)/B)$ does not usually give an approximately
trivial extension.
When $A$ is stably finite, both $K_0(A)$ and $K_0(M(B))$
have nice order while $K_0(M(B)/B)$ has no useful order. Even if we know
which \hm\, $\beta: K_0(A)\to K_0(M(B)/B)$ can be lifted to a \hm\,
from $K_0(A)$ to $K_0(M(B)),$ the lifting may not be positive.
In this paper, at least for some special cases,
we give a precise condition for an element in $KL(A, M(B)/B)$
to be represented by approximately trivial extensions.

This paper is organized as follows.

Section 1: Preliminaries\\
This section is a preparation for the rest of the paper which contains
a computation of $K$-theory for $M(B)$ and $M(B)/B$ for
$\sigma$-unital simple \CA\, with real rank zero, stable rank one,
weakly unperforated $K_0(B)$ and with a continuous scale.
We also point out that  $M(B)/B$ is simple (and purely infinite)
if and only if $B$ has a continuous scale (if $B\not\cong {\cal K}$).

Section 2: Monomorphisms from $A\otimes {\cal O}_2$ into a purely
infinite simple \CA\\
This  section studies \hm s from $A\otimes {\cal O}_2$ into a purely infinite
simple \CA.

Section 3: Approximately unitarily equivalent extensions\\
We show that, if $B$ is a non-unital and
$\sigma$-unital simple \CA\, with a continuous scale,
two monomorphisms from $A$ to
$M(B)/B$ are approximately unitarily equivalent if and only if
they induce the same element in $KL(A, M(B)/B).$

Section 4: ${\bf Ext}_{ap}(A,B)$\\
In this section, under the assumption that $A$ satisfies the UCT, we give a bijection $\Gamma: {\bf Ext}_{ap}(A,B)\to KL(A, M(B)/B).$

Section 5: Examples\\
In this section, we present a few examples which show that the bijection $\Gamma$
may not answer all questions about these extensions. For example, we show that
the zero element in $KL(A,M(B)/B)$ does not represent an approximately trivial extension in general.

Section 6: Quasidiagonal extensions --- general and infinite cases\\
This section discusses quasidiagonal extensions. Without assuming the UCT, we give
a general $K$-theoretical necessary condition for an essential extension
to be quasidiagonal.
We also show that for any separable exact \CA\, $A,$ there exist quasidiagonal extensions of $A$ by any $\sigma$-unital purely infinite simple \CA s.

Section 7: Quasidiagonal extensions --- finite case\\
Let $A$ be a separable amenable \CA\,
and $B$ be a $\sigma$-unital \CA\, admitting an approximate identity consisting of projections and having the property (SP). We show that, if $A$ is a quasidiagonal \CA, then there exists an essential quasidiagonal extension.
If, in addition that $B$ is also assumed to be a quasidiagonal \CA, then
the condition that $A$ is quasidiagonal is also necessary. When $B$ is a $\sigma$-unital simple \CA\, with real rank zero, stable rank one, weakly unperforated $K_0(B)$ and with a continuous scale, we present a $K$-theoretical
necessary and sufficient condition for an essential extension to be quasidiagonal for a class of separable quasidiagonal amenable \CA s.

Section 8: Approximately trivial extensions\\
In the last section, we give a general $K$-theoretical necessary condition
for essential extensions to be approximately trivial. Combining this condition
with the results in section 7, we show that there are essential quasidiagonal extensions that are {\it not} approximately trivial extensions. We also show
how to use the bijection $\Gamma$ to determine which essential extensions are
approximately trivial at least in some special cases.

{\bf Acknowledgement} Part of this work was done when
the author was visiting East China Normal University.
This work was partially supported by a grant from National
Science Foundation of U. S. A.

\section{Preliminaries}

Throughout this paper, we will use the following conventions:

1) An ideal of a \CA\, is always a closed two-sided ideal.

2) By a unital \SCA\, $C$ of a unital \CA\, $A$ we mean $C\subset
A$ and $1_C=1_A.$

3) If $p$ and $q$ are two projections in a \CA\, $A,$ we say $p$ is
equivalent to $q$ if there exists a partial isometry $v\in A$ such
that $v^*v=p$ and $vv^*=q.$

4) Let $A$ and $B$ be two \CA s and $L_1, L_2: A\to B$ be two
maps. Let $\ep>0$ and ${\cal F}\subset A$ be a subset. We write
$$
L_1\approx_{\ep} L_2\,\,\,\, {\rm on}\,\,\, {\cal F}
$$
if
$$
\|L_1(a)-L_2(b)\|<\ep\,\,\,\,{\rm for \,\,\, all}\,\,\, a\in {\cal
F}.
$$

Suppose that $A$ and $B$ are unital
and $L_1(1)$ and $L_2(1)$ are projections.
If there is an isometry $s\in B$ such that $s^*L_2(1)s=L_1(1),$ $sL_1(1)s^*=L_2(1)$ and
$$
{\rm ad}\, s\circ L_2\approx_{\ep} L_1\,\,\,\, {\rm on}\,\,\, {\cal F},
$$
we will write
$$
L_2 {\sim}_{\ep} L_1\,\,\,\, {\rm on}\,\,\, {\cal F}.
$$

5) A separable \CA\, $A$ is said to be ${\it amenable}$ (or
nuclear), if for any $\ep>0$ and finite subset ${\cal F}\subset
A,$ there exists a finite dimensional \CA\, $C$ and two \morp s
$L_1: A\to C$ and $L_2: C\to A$ such that
$$
L_2\circ L_1\approx_{\ep} \id_A \,\,\,{\rm on}\,\,\, {\cal F}.
$$

6) Let $A$ be a separable amenable \CA. We say $A$ satisfies the
Universal Coefficient Theorem (UCT) and write $A\in {\cal N},$ if
for any $\sigma$-unital \CA\, $C,$ the map $\gamma: KK(A,C)\to
{\rm Hom}(K_*(A), K_*(C))$ is surjective and the map $\kappa: {\rm
ker}\gamma\to ext(K_*(A), K_{*}(C))$ is an isomorphism, i.e., there
is a short exact sequence:
$$
0\to ext(K_*(A), K_*(C)) \to KK(A,C)
{\stackrel{\gamma}{\to}}{\rm Hom}(K_*(A), K_*(C))\to 0.
$$
 If $h: A\to C$ is a \hm\, then $h$ gives an element $[h]$ in $KK(A,C).$

7) An extension $0\to B\to E\to A\to 0$ of \CA s is said
to be essential if
$$
\{e\in E: eb=be=0 \,\,\,{\rm for\,\,\, all}\,\,\, b\in B\}=\{0\}.
$$
If $E$ is an essential extension of $A$ by $B$ as above, then
it is determined by a monomorphism $\tau: A\to M(B)/B$ and
$E=\pi^{-1}\circ \tau(A),$ where $\pi: M(B)\to M(B)/B$ is
the quotient map.

\begin{NN}
\end{NN}
\vspace{-0.1in}
We start with the following essential extensions:
$$
0\to A\to E\to C(D)\to 0,
$$
where  $D$ is the unit disk and $A=B\otimes {\cal K}$ and where
$B$ is a separable unital simple AF-algebra with a unique tracial
state. For example $B$ may be  a UHF-algebra. Let $I$ be the unique
proper  ideal of $M(A)$
 which contains $A$ (3.2 in \cite{Ell1}). Denote $J=\pi(I),$ where
$\pi: M(A)\to M(A)/A$ is the quotient map. Let $p\in
M(A)/A\setminus J$ be a projection such that $1-p\in J$ is a
non-zero projection. To see such a projection exists, one takes a
projection $q\in M(A)\setminus A$ with finite trace. Then $q\in
I.$ Define $p=1-\pi(q).$ It follows from 1.17 (4) in \cite{Ln5}
that $K_1(M(A)/I))={\mathbb R}.$ It is known that $M(A)/I$ is
purely infinite and simple (see \cite{Zh2}). Thus there is a
unitary $u\in p(M(A)/A)p$ such that ${\bar \pi}(u)$ is not in
$U_0(M(A)/I),$ where ${\bar \pi}: M(A)/A\to M(A)/I$ is the
quotient map. Let $y\in (1-p)M(A)/A(1-p)$ with the spectrum
$sp(y)=D.$ Set $x=u+y.$ Define $\tau: C(D)\to M(A)/A$ by
$\tau(f)=f(x)$ for $f\in C(D).$ It is easy to see that $\tau$ is
not trivial nor it is approximately trivial. However, it is
known that  $Ext(C(D),A)=KK^1(C(D), A)=\{0\}.$ So certainly in this case
$KK^1(C(D), A)$ can not be used to understand extensions of $C(D)$
by $A.$  Clearly the complicity of the extension is caused by the
fact that $M(A)/A$ is not simple. One can easily imagine that when
the ideal structure of $M(A)/A$ is more complicated, equivalent
classes of extensions will be hard to compute if it is even
possible to compute. The success of the BDF-theorem for the
classification of extensions by ${\cal K}$ depends on the fact
that $M({\cal K})/{\cal K},$ the Calkin algebra, is simple. In this
paper, we will therefore consider only those essential extensions
by a simple \CA\, $A$ such that $M(A)/A$ is simple.

So the question is:

\begin{NN}
When is $M(B)/B$ simple?
\end{NN}
\vspace{-0.1in}
 Let $B$ be a $\sigma$-unital simple \CA .  Recall (\cite{Ln1})
that $B$ is said to have a {\it continuous scale} if for any
approximate identity $\{e_n\}$ of $B$ satisfying
$e_{n+1}e_n=e_ne_{n+1}=e_n$  and any nonzero positive element $a\in
B_+,$ there exists an integer $N>0$ such that
$$
(e_m-e_n) \lesssim  a,\,\,\, {\rm for}\,\,\, m>n\ge N
$$
i.e., there exists a sequence of elements $r_n\in B$ such that
$r_k^*ar_k\to e_m-e_n$ for all $m>n\ge N.$ It should be noted that
if $p$ and $q$ are projections and $p\lesssim q,$ then $p$ is
equivalent to a projection $q'\le q.$

It is proved in \cite{Ln1} that, for non-unital separable simple
\CA\, $B\not\cong {\cal K},$  $M(B)/B$ is simple if $B$ has
a continuous scale.  Recently we have proven the following:

\vspace{0.1in}

\begin{thm}{\rm (\cite{Lnnew})}\label{0Tadd}
Let $A\not\cong {\cal K}$ be a non-unital and $\sigma$-unital simple
\CA. The following are equivalent:

{\rm (1)} $A$ has a continuous scale;

{\rm (2)} $M(A)/A$ is simple,

{\rm (3)} $M(A)/A$ is a purely infinite simple \CA.

\end{thm}

Clearly  every (non-unital) $\sigma$-unital purely infinite simple
\CA\, has a continuous scale. Essential extensions of separable \CA
s $A$ which satisfy the UCT by a non-unital separable purely
infinite simple \CA\, $B$ is classified by  $KK^1(A, B)$ by
Kirchberg's absorbing theorem (\cite{K1}).

In this  paper we will focus on essential extensions by a
$\sigma$-unital simple \CA\, with real rank zero, stable rank one,
weakly unperforated $K_0$ and a continuous scale.


Suppose that $B$ is a non-unital separable simple \CA\, with real
rank zero, stable rank one and weakly unperforated $K_0(B).$  Fix
any nonzero projection $e\in B.$ Denote by $T$ the set of those
quasi-traces $\tau$ on $B$ such that $\tau(e)=1.$ Note that $T$ is (weak
$*$-) compact convex set.  Let $a\in M(B)_+.$ Define
$\hat{a}(\tau)=\tau(a)$ for $\tau\in T.$ Then $\hat{a}$ is a lower
semi-continuous  affine function on $T.$ If $a\in A,$ then
$\hat{a}$ is continuous.

To see examples of simple \CA s with continuous scale, we quote
the following result (\cite{Lnnew}). It also justifies the
terminology ``continuous scale".

\begin{thm}\label{0Tcs}
Let $A$ be a non-unital but  $\sigma$-unital simple \CA\, with
real rank zero, stable rank one and weakly unperforated $K_0(A).$
Let $1$ be the identity of $M(A).$ Then $A$ has a continuous scale
if and only if $\hat{1}(\tau)=\tau(1)$ for $\tau\in T$ is a
continuous function on $T.$
\end{thm}

It is also proved in \cite{Lnnew} that given any  separable simple
\CA\, $A,$ there is a non-unital hereditary \SCA\, $B$ such that
$B$ has a continuous scale. In particular, $M(B)/B$ is a purely
infinite simple \CA. Note that $B\otimes {\cal K}\cong A \otimes
{\cal K}$ and $B$ may not have any non-trivial projections.
Furthermore, $B$ may contain both infinite
 and finite projections (given by R\o rdam (\cite{Ro5}).

\begin{Def}\label{0Daff}
{\rm Let $T$ be a compact convex set. A function $f$ defined on
$T$ is
said to be affine if \\
$f(a\xi+(1-a)\zeta)=af(\xi)+(1-a)f(\zeta)$
for all $\xi,\zeta\in T$ and $0\le a\le 1.$  We denote
$$
\mathit{Aff}(T)=\{f\in C(T): f\,\,\,{\rm is\,\,\, affine}\}.
$$
If $f, g\in \mathit{Aff}(T)$ and $f(t)>g(t)$ for all $t\in T,$ we will
write $f\gg g.$ Denote
$$
 Aff(T)_{++}=\{f\in Aff(T): f\gg 0, {\rm or}\,\,\,
f=0\}. $$

Let $B$ be a simple \CA\, with real rank zero, stable rank one and
weakly unperforated $K_0(B).$ Fix a nonzero projection $e\in B.$
Denote by $T$ the set of those  traces $\tau$  defined on $B$ such
that $\tau(e)=1.$ Define
$$
\rho_B: K_0(B)\to Aff(T)
$$
by $\rho_B([p])(\tau)=\tau(p)$ for  projections $p\in M_m(B),$
$m=1,2,...$ It is known that $\rho_B$ is a positive \hm\, from
$(K_0(B), K_0(B)_+)$ to $(Aff(T),Aff(T)_{++})$ (see \cite{BH}). In
fact (by \cite{BH}), $[p]\ge [q]$ and $[p]\not=[q]$ if and only if
$\tau(p)>\tau(q)$ for all $\tau\in T.$

}
\end{Def}

The following was first proved in \cite{Ln2} in 1991.

\begin{thm}\label{0T1}
Let $e\in B$ be a nonzero projection. Let
$$
T=\{ \tau: \tau(e)=1, \tau\,\,\, {\rm is\,\,\, a\,\,\,
trace\,\,\, defined \,\,\,
on}\,\,\, B\}.
$$
 Then

{\rm (1)} $(K_0(M(B)), K_0(M(B))_+)=(Aff(T), Aff(T)_{++});$

{\rm (2)} two projections $p$ and $q$ in $M(B)\setminus B$ are
equivalent if $\tau(p)=\tau(q)$ for all $\tau\in T;$

{\rm (3)} for any $f\in Aff(T)_{++},$ there is a projection
$p\in M_k(M(B))\setminus M_k(B)$ (for some $k\ge 1$) such that
$\hat{p}=f;$ and

{\rm (4)} $K_1(M(B))=\{0\}$ and $U(M(B))=U_0(M(B)).$

\end{thm}

\begin{proof}
  Since $B$ has real rank zero, we obtain an approximate identity $\{e_n\}$ for
$B$ consisting of projections (with $e_0=0$).  Let $p\in M(B)$ be
a projection. It follows from Theorem 4.1 in \cite{Zh1} that we
may assume that $p=\sum_{n=1}^{\infty}p_n,$ where $p_n\le
e_{n+1}-e_{n}.$ Since $\hat{1}$ is continuous on $T,$ by the Dini
Theorem, $\rho_B(e_n)$ converges to
$\hat{1}$ uniformly on $T.$ This implies that $\hat{p}$ is also
continuous. Define $\rho: K_0(M(B))\to Aff(T)$ by defining
$\rho([p])=\hat{p}.$ It is clear that $\rho$ is a well-defined \hm.

We now prove (2). It follows from Theorem 4.1 in \cite{Zh1} again that we may assume that $p=\sum_{n=1}^{\infty}p_n$ and
$q=\sum_{n=1}^{\infty}q_n,$ where $p_n, q_n\le e_{n+1}-e_{n}$ and
the sum converges in the strict topology.
 Without loss of generality, we may assume that $p_n$ and $q_n$ are not zero.
Since $B$ is simple, we have $\widehat{p_1}\ll \hat{p}=\hat{q}.$
Since $\sum_{n=1}^{\infty}\hat{q_n}$ converges uniformly on $T,$
there is $n_1>1$ such that
$\sum_{j=1}^{n_1}\widehat{q_j}\gg\widehat{p_1}.$ It follows
from III2.2 and III2.3  in \cite{BH} that there is a partial isometry $v_1\in B$
such that
$$
v_1^*v_1=p_1\andeqn v_1v_1^*\le \sum_{j=1}^{n_1}q_j.
$$
There is $m_1>1$ such that
$$
\tau(\sum_{j=2}^{m_1}p_j)>\tau(\sum_{j=1}^{n_1}q_j -v_1v_1^*)\,\,\,\,
{\rm for\,\,\,all}\,\,\,\tau\in T.
$$
It follows that there is a partial isometry $u_1\in B$ such that
$$
u_1^*u_1=\sum_{j=1}^{n_1}q_j -v_1v_1^* \andeqn u_1u_1^*\le
\sum_{j=2}^{m_1}p_j.
$$
Put $w_1=v_1+u_1^*.$ Then
$$
\sum_{j=1}^{m_1}p_j >w_1^*w_1\ge p_1\andeqn \sum_{j=1}^{n_1}q_j=w_1w_1^*\ge q_1.
$$
By induction one constructs a sequence of partial isometries
$w_k\in B$ ($w_0=0$) such that
$$
\sum_{j=1}^{m_k}p_j-\sum_{j=1}^{k-1}w_j^*w_j >w_k^*w_k\ge
\sum_{j=1}^{m_{k-1}}p_j-\sum_{j=1}^{k-1}w_j^*w_j \andeqn
\sum_{j=1}^{n_k}q_j-\sum_{j=1}^{k-1}w_jw_j^*=w_kw_k^*\ge
\sum_{j=1}^{n_{k-1}}q_j -\sum_{j=1}^{k-1}w_jw_j^*,
$$
where $\{n_k\}$ and $\{m_k\}$ are increasing sequences of positive integers.

Set $W=\sum_{k=1}^{\infty}w_k.$ One checks that the sum
converges in the strict topology and $W$ is a partial isometry in
$M(B).$ One then verifies that
$$
W^*W=p\andeqn WW^*=q.
$$

This proves (2).

Note that (2) also implies that $\rho$ is injective. In fact, if
$p\in M(B)\setminus B$ and $q\in B$ are two projections such that
$\tau(p)=\tau(q)$ for all $\tau\in T,$ then $p\oplus 1$ and
$q\oplus 1$ are both in $M(B)\setminus B.$ So, by (2), $p\oplus 1$
and $q\oplus 1$ are equivalent. Therefore $\rho$ is injective.

To see $\rho$ is surjective, let $f\in Aff(T).$ We need to show
that $f$ is in the image of $\rho.$ It is clear that it suffices
to prove the case for $f\gg 0.$ So we may assume that $f\gg 0.$ We
claim that there exists a sequence of positive functions
$\{f_n\}$ in $\rho_B(K_0(B))$ such that $f_n\ll f_{n+1}$ such that
$f_n\to f$ uniformly on $T.$

Let $d_0=\inf\{f(t): t\in T\}.$ Then $d_0>0.$

Let $d_0/2>\ep>0.$ Since $\rho_A(K_0(B))$ is dense in $Aff(T)$
(see \cite{BH}), there is $g_1\in \rho_B(K_0(B)_+)$ such that
$$
\|g_1-(f-\ep)\|<\ep/4.
$$
Therefore $(1-\ep/2)f\ll g_1\ll f.$ By applying the same argument to
the function $f-g_1,$ we obtain $g_2\in \rho_B(K_0(B)_+)$ such
that $(1-\ep/4)(f-g_1)\ll g_2\ll f-g_1.$ Note that $g_1+g_2\in
\rho_B(K_0(B)_+).$ From this the claim follows.

Now we will show that $f$ is in the image of $\rho.$ By replacing
$f$ by $f-f_n$ if necessary, we may assume that $f\ll
\widehat{e_1}.$ There are projections $r_n\in B$ such that
$\rho(r_n)=f_n-f_{n-1}$ (with $f_0=0$). Since $f\ll \hat{e_1}$ and
$f_n$ converges to $f$ uniformly on $T,$ we may assume that
$$
\sum_{k=1}^{n_1}(f_k-f_{k-1})\ll \rho(e_1)\andeqn
\sum_{k=n_1+1}^{\infty} (f_k-f_{k-1})\ll \rho(e_2-e_1).
$$
Thus  we may assume that
$$
\sum_{k=1}^{n_1}r_k\le e_1\andeqn
\sum_{k=n_1+1}^m r_k\le e_2-e_1
$$
for all $m>n_1.$
We obtain $n_2>n_1$ such that
$$
\sum_{k=n_2}^{\infty} (f_k-f_{k-1})\ll \rho(e_3-e_2)
$$
Therefore we also assume that
$$
\sum_{k=n_1+1}^{n_2} r_k\le e_2-e_2\andeqn
\sum_{k=n_2+1}^mr_k\le e_3-e_2
$$
for all $m>n_2.$
By induction, we obtain an increasing sequence of integers $\{n_k\}$ such that
$$
\sum_{i=n_k+1}^{n_{k+1}}r_i\le e_{k+1}-e_k,\,\,\, k=1,2,....
$$
Note that this implies that $\sum_{n=1}^{\infty} r_n$ converges in
the strict topology to a projection $p\in M(B).$ It is easy to see
that $\hat{p}=f.$ This proves (1) as well as (3).

Finally we note  $pBp$ has real rank zero for all
projection $p\in M(B)$ and by \cite{Lnep}
$cer(pAp)\le \pi.$ Thus  (4)
has been proved in Lemma 3.3 in \cite{Lngw}.

\end{proof}

 \begin{Cor}\label{0T2}
Let $B$ be as in \ref{0T1}. Then

{\rm (1)} $K_1(M(B)/B)={\rm ker}\rho_B,$ and

{\rm (2)} there is a short exact sequence
$$
0\to Aff(T)/\rho_A(K_0(B))\to K_0(M(B)/B)\to K_1(B)\to 0.
$$
\end{Cor}

\begin{proof}

From the following six-term exact sequence
$$
\begin{array}{ccccc}
K_0(B) & {\stackrel{\rho_B}{\to}} &K_0(M(B)) & \to & K_0(M(B)/B)\\
\uparrow && &&\downarrow\\
K_1(M(B)/B)& \leftarrow & K_1(M(B)) & \leftarrow & K_1(B)\\
\end{array}
$$
we obtain, by \ref{0T1},
$$
\begin{array}{ccccc}
K_0(B)& {\stackrel{\rho_B}{\to}} & Aff(T)& \to & K_0(M(B)/B)\\
\uparrow && &&\downarrow\\
K_1(M(B)/B)& \leftarrow & 0 & \leftarrow & K_1(B)\\
\end{array}
$$
This six-term exact sequence unsplices into
$$
K_1(M(B)/B)={\rm ker}\rho_B \andeqn 0\to Aff(T)/\rho_B(K_0(B))\to
K_0(M(B)/B)\to K_1(B)\to 0.
$$

\end{proof}

\begin{Remark}
{\rm It should be noted that $Aff(T)$ as an ordered group does not depend on the choice
of the non-zero projection $e.$
In what follows when we write $Aff(T)$ it is understood that
the projection $e$ is fixed.}
\end{Remark}

The following fact will be used in this paper.

\begin{Prop}\label{0Pph}
Let $G$ be a dense ordered subgroup of ${\mathbb R}$ containing $1$
and let $T$ be a Choquet simplex.
Suppose that $h: G\to Aff(T)$ is a positive \hm\, with
$h(1)=a.$ Then
$$
h(z)=za\,\,\, for\,\,\, all\,\,\, z\in G.
$$
\end{Prop}

\begin{proof}
Since $G$ is dense in ${\mathbb R},$ there exists a
sequence $g_n>0$ in $G$ such that $g_n\to 0.$
We may assume that $ng_n<1$ for all $n.$
Therefore
$$
nh(g_n)\le a\,\,\, n=1,2,....
$$
It follows that $h(g_n)\to 0.$
Let $f_n>0$ in $G$ such that $f_n\to 0.$
Thus, for each $m,$ there exists $N(m)$ such that
$f_n\le g_m$ if $n\ge N(m).$ This implies
that $h(f_n)\to 0.$
Thus $h$ is continuous.
For each nonzero integer $m,$ define ${\tilde h} (1/m)=a/m.$
Then one checks that ${\tilde h}$ is a positive \hm\,
from ${\mathbb Q}G$ to $Aff(T).$
The same argument above shows that ${\tilde h}$ is also
continuous. Fix $z\in G.$ Suppose that $r_n\in {\mathbb Q}$ such
that $r_n\to z.$ Then ${\tilde h}(r_n)\to h(z),$
or $r_na\to h(z).$  Therefore $h(z)=az.$
\end{proof}

The following example shows
that even in the case that $M(B)/B$ is simple, in general,
$KK^1(A,B)$ can not be used to give a meaningful
description of extensions of $A$ by $B.$

\begin{Ex}\label{0E2}
{\rm Let $A$ be a unital simple AF-algebra and $B$ be a
$\sigma$-unital simple \CA\, with real rank zero, stable rank one,
weakly unperforated $K_0(B)$ and a continuous scale. Let $r$ be an
irrational number and $D_r=\{m+nr: m,n\in {\mathbb Z}\} $
($={\mathbb Z}\oplus {\mathbb Z}r$). Suppose $K_0(A)=D_r\oplus
{\mathbb Z}.$ Define
$$
K_0(A)_+=\{ x+z: x\in D_r, x>0\}\cup\{(0,0)\}.
$$
Suppose that there is a group \hm\, $\theta: K_0(A)\to
Aff(T)/\rho_B(K_0(B))$ such that $\theta((0,1))\not=0,$ where
$(0,1)$ denotes the generator of the last summand ${\mathbb Z}$ of
$K_0(A).$  This gives a group \hm\, $\alpha: K_0(A)\to
K_0(M(B)/B)$ which maps $(0,1)$ to $r.$ Let $\Phi: Aff(T)\to
Aff(T)/\rho_B(K_0(B))$ be the quotient map. Then $\theta=\Phi\circ
\alpha$ gives one such \hm. Since $M(B)/B$ is a purely infinite
simple \CA, it is easy to construct a \hm\, $\tau: A\to M(B)/B$
such that $\tau_{*0}=\theta$ (see for example \ref{IITM}). This
$\tau$ gives an essential extension of $A$ by $B$ which gives an
element in $KK^1(A,B).$ Let $E$ be the \CA\, determined by $\tau$
we have the following commutative diagram:
$$
{\small \put(-160,0){$K_0(B)$} \put(20,0){$K_0(E)$}
\put(180,0){$K_0(A)$} \put(-85,-40){$K_0(B)$}
\put(0,-40){$K_0(M(B))$} \put(105,-40){$K_0(M(B)/B))$} \put(-85,
-70){$K_1(M(B)/B)$} \put(27,-70){$0$} \put(115,-70){$K_1(B)$}
\put(-160,-110){$K_1(A)$} \put(15,-110){$K_1(E)$}
\put(180,-110){$K_1(B)$} \put(-120, 2){\vector(1,0){95}}
\put(70,1){\vector(1,0){95}} \put(-123,-3){\vector(1,-1){30}}
\put(30,-3){\vector(0,-1){25}} \put(180,-2){\vector(-1,-1){30}}
\put(-45,-38){\vector(1,0){35}} \put(70,-38){\vector(1,0){25}}
\put(-147, -90){\vector(0,1){85}} \put(-75,-60){\vector(0,1){15}}
\put(125, -45){\vector(0,-1){15}} \put(190,-7){\vector(0,-1){85}}
\put(15,-68){\vector(-1,0){30}} \put(95,-68){\vector(-1,0){25}}
\put(-123,-102){\vector(1,1){30}} \put(175,
-105){\vector(-1,1){30}} \put(30,-100){\vector(0,1){25}} \put(-5,
-108){\vector(-1,0){100}} \put(170,-108){\vector(-1,0){95}}
\put(-111,-13){$=$}  \put(151, -13){$\tau_{*0}$} \put(160,
-88){$=$} }
$$

Since the image of $\tau_{*0}$ is in $Aff(T)/\rho_B(K_0(B)),$
by (2) in  \ref{0T2} and from
the above diagram one concludes that the map from $K_0(A)$ to
$K_1(B)$ is zero. Since $K_1(A)=\{0\},$ one further concludes that
the map from $K_1(A)$ to $K_0(B)$ is also zero. By the Universal
Coefficient Theorem, one computes that $[\tau]\in ext_{\mathbb
Z}(K_0(A), K_0(B)).$ Using the map $\alpha$ or using the fact that
$K_0(A)$ is finitely generated free group, $[\tau]=0$ in
$KK^1(A,B).$ However, there is no \hm\, $h: A\to M(B)$ such that
$(\pi\circ h)_{*0}=\theta.$ Otherwise, since $h_{*0}$ is positive,
it maps ${\rm ker}\rho_A={\mathbb Z}$ into zero. It would imply
that $\tau_{*0}$ maps ${\rm ker}\rho_A$ to zero. But we
constructed otherwise. Therefore $\tau$ is not trivial.
Furthermore it can not be even approximately trivial (see \ref{IITapptr}
below). This shows that even in the case that $M(B)/B$ is simple,
$KK^1(A,B)$ can not be used to give a good description of extensions
of $A$ by $B.$ }
\end{Ex}

\section{Monomorphisms from $A\otimes {\cal O}_2$ into a purely infinite
simple \CA}

\begin{Def}\label{IIDfilter}
{\rm Recall that a family ${\omega}$  of subsets of ${\mathbb N}$
is an {\it ultrafilter} if

(i) $X_1,...,X_n\in \omega$ implies $\cap_{i=1}^n X_i\in \omega,$

(ii) $\O\not\in \omega,$

(iii) if $X\in\omega$ and $X\subset Y,$ then $Y\in\omega$ and

(iv) if $X\subset {\mathbb N}$ then either $X$ or ${\mathbb
N}\setminus X$ is in $\omega.$

An ultrafilter is said to be {\it free}, if $\cap_{X\in \omega} X
=\emptyset.$ The set of free ultrafilters is identified with elements in
$\beta {\mathbb N}\setminus {\mathbb N},$ where $\beta{\mathbb N}$
is the Stone-Cech compactification of ${\mathbb N}.$

 A sequence $\{x_n\}$  (in a normed space for example) is
said to converge to $x_0$ along $\omega,$ written
$\lim_{\omega}x_n=x_0,$ if for any $\ep>0$ there exists
$X\in\omega$ such that $\|x_n-x_0\|<\ep$ for all $n\in X.$

Let $\{B_n\}$ be a sequence of \CA s. We write
$l^{\infty}(\{B_n\})$ for the \CA\, $\prod_{n=1}^{\infty} B_n.$
Fix an ultrafilter $\omega.$ The ideal of $l^{\infty}(\{B_n\})$
which consists of those sequences $\{a_n\}$ in $l^{\infty}(\{B_n\})$
such that $\lim_{\omega}\|a_n\|=0$ is denoted by $c_{\omega}(\{B_n\}).$
Define
$$
q_{\omega}(\{A_n\})=l^{\infty}(\{B_n\})/c_{\omega}(\{B_n\}).
$$
If $B_n=B, n=1,2,...,$ we use $l^{\infty}(B)$ for
$l^{\infty}(\{B_n\}),$ $c_{\omega}(B)$ for $c_{\omega}(\{B_n\})$
and $q_{\omega}(A)$ for $q_{\omega}(\{A_n\}),$ respectively. }
\end{Def}

\begin{Lem}\label{NIL1}
Let $A$ be a separable \CA\, and $\{B_n\}$ be a sequence of unital \CA s.
Let $\omega\in \beta{\mathbb N}\setminus {\mathbb N}.$ Suppose that
$\psi_m, \phi_m: A\to B_m$  are two bounded sequences of maps such that
$\psi=\pi\circ \{\psi_m\}, \phi=\pi\circ \{\phi_m\}: A\to q_{\omega}(\{B_n\})$
are two \hm s, where
$\pi: l^{\infty}(\{B_m\})\to q_{\omega}(\{B_n\})$ is the quotient map.

{\rm (1)} Suppose that there are isometries $u_n\in q_{\omega}(\{B_n\})$ such that
$$
\lim_{n\to\infty}\|u_n^*\psi(a) u_n -\phi(a)\|=0\,\,\,{\rm for \,\,\, all}\,\,\,
a\in A.
$$
Then there is an isometry $w\in q_{\omega}(\{B_n\})$ such that
$$
w^*\psi(a)w=\phi(a)\,\,\,{\rm for \,\,\, all}\,\,\,
a\in A.
$$

{\rm (2)} Suppose that
$\psi$ and $\phi$ are  approximately
unitarily equivalent in $q_{\omega}(\{B_n\}).$ Then
they are unitarily equivalent.
\end{Lem}
\begin{proof}
Suppose that there is a sequence of isometries $u_n\in
q_{\omega}(\{B_n\})$ such that
$$
\lim_{n\to\infty}\|u_n^*\psi(a)u_n-\phi(a)\|=0
$$
for all $a\in A.$
Let $\{a_n\}$ be a dense sequence of $A.$ By passing to a subsequence if necessary, we may assume that
$$
\|u_n^*\psi(a_j)u_n-\phi(a_j)\|<1/n,\,\,\,j=1,2,...,n.
$$
It follows from 6.2.4 in \cite{Rob} that there exists, for each
$n,$ a sequence of isometries $u_m^{(n)}\in B_n$ such that $\pi(\{u_m^{(n)}\}=u_n,$ where $\pi: l^{\infty}(\{B_n\})\to q_{\omega}(\{B_n\})$
is the quotient map.
For each $n,$ there exists $X_n\in \omega$ such that for $m\in X_n,$
$$
\|(u_m^{(n)})^*\psi_m(a_j)u_m^{(n)}-\phi_m(a_j)\|\le 1/n,j=1,2,...,n.
$$
Since $\omega$ is free, there is for each $j,$ $Y_j'\in \omega$ such that $j\not\in Y_j'.$ Let $Y_j=\cap_{1\le k\le j} Y_j'.$ Then $Y_j\in \omega$
and $\{1,2,...,j\}\cap Y_j=\emptyset.$
Let $Z_k'=X_k\cap Y_k.$ Then $Z_k'\in \omega$ and $\cap_{k\ge N}Z_k'=\emptyset,$
$N=1,2,....$ Put $Z_k=\cap_{1\le j\le k}Z_j'.$ Then $Z_k\in \omega$ and
$Z_1\supset Z_2\supset \cdots \supset Z_k\supset \cdots.$ Moreover, $Z_k\subset X_k,$ $k=1,2,....$

Define $l(m)$ as follows. If $m\in Z_k\setminus Z_{k+1},$ define $l(m)=k,$ $k=1,2,...;$ and if $m\not\in Z_1,$ define $l(m)=m,$ $m=1,2,....$
Put $w_m=u_m^{(l(m))}\in B_m$ and $w=\pi(\{w_m\}).$
Then, for any $\ep>0$ and $j,$ let $k>0$ be an integer such that $1/k<\ep$ and $j\le k.$
If $m\in Z_k,$ then $m\in Z_{k'}\setminus Z_{k'+1}$ for some
$k'\ge k.$ Thus $w_m=u_m^{(k')}$ and $m\in X_{k'}.$
Therefore
$$
\|w_m^*\psi_m(a_j)w_m-\phi_m(a_j)\|=\|u_m^{(k')}\psi_m(a_j)u_m^{(k')}-\phi_m(a_j)\|<1/k'\le 1/k<\ep
$$
for all $j=1,2,...,k.$
This implies that
$$
\lim_{\omega}\|w_m^*\psi_m(a_j)w_m-\phi_m(a_j)\|=0
$$
for all $j.$
Hence
$$
w^*\psi(a_j)w=\phi(a_j)\,\,\,j=1,2,....
$$
Since $\{a_n\}$ is dense in $A,$
$$
w^*\psi(a)w=\phi(a)\,\,\,{\rm for \,\,\, all}\,\,\, a\in A.
$$
This proves (1).

To prove (2), we note that if $u_n$ are unitaries, so
is $w.$
\end{proof}

\begin{Lem}{\rm (\cite{KP}, Proposition 1.4)}\label{NIL3}
Let $A$ be a unital separable \SCA\, of a unital purely infinite
simple \CA\, $B$ and let $\phi: A\to B$ be an amenable \morp. Then
for any finite subset ${\cal F}\subset A$ and any $\ep>0$ there is
a non-unitary isometry $s\in B$ such that
$$
\|s^*as-\phi(a)\|<\ep\,\,\,\,for\,\,\, all\,\,\, a\in {\cal F}.
$$
\end{Lem}

\begin{proof}
By the assumption that $\phi$ is amenable there
are \morp s $\sigma: A\to M_n({\mathbb C})$ (for some integer $n>0$)
and $\eta: M_n({\mathbb C})\to B$ such that
$$
\phi\approx_{\ep} \eta\circ \sigma\,\,\,\,{\rm on}\,\,\, {\cal F}.
$$
We may therefore assume that $\phi=\eta\circ \sigma.$
It is well known (see for example, 2.3.5 in \cite{Lnb})
 that there exists a \morp\,
${\tilde \sigma}: B\to M_n({\mathbb C})$ such that
${\tilde \sigma}|_{A}=\sigma.$
Define ${\tilde \phi}=\eta\circ {\tilde \sigma}.$
Since now $B$ is purely infinite and ${\tilde \phi}$ is amenable,
the lemma follows immediately from 1.4 in \cite{KP}.
\end{proof}

\begin{Cor}\label{NIL4}
Let $B$ be a unital purely infinite simple \CA\,
and $A$ be a separable amenable
\CA. Let $\phi,\psi: A\to B$ be two monomorphisms. Then
there are two sequences of isometries $s_n$ and $w_n$ in $B$ such that
$$
\lim_{n\to\infty}\|s_n^*\phi(a)s_n-\psi(a)\|=0 \,\,\, and \,\,\,
\lim_{n\to\infty}\|w_n^*\psi(a)w_n-\phi(a)\|=0
$$
for all $a\in A.$
\end{Cor}

The proof of the following proposition is exactly the same as
that of 6.2.6 in \cite{Rob}.

\begin{Prop}\label{NIPps}
Let $B_n$ be a sequence of purely infinite simple \CA s.
Then $q_{\omega}(\{B_n\})$ is a purely infinite simple \CA\, for every free ultrafilter $\omega.$
\end{Prop}

\begin{Prop}{\rm (3.4 in \cite{KP})}\label{NIP1}
Let $B_n$ be a sequence of unital purely infinite simple \CA s and $A$ be a
unital separable amenable simple \CA. Suppose
that $j: A\to q_{\omega}(\{B_n\})$ is a monomorphism.
Then the relative commutant $j(A)'$ in $q_{\omega}(\{B_n\})$
is a unital purely infinite simple \CA.
\end{Prop}

\begin{proof}
Let $a\in j(A)'$ be a nonzero positive element with $\|a\|=1.$ It
suffices to show that there is an isometry $s\in j(A)'$ such that
$s^*as=1.$ Let $X=sp(a)\subset [0,1]$ and define two  \hm s
$\phi,\psi: C(X)\otimes A\to q_{\omega}(B)$ by
$$
\phi(f\otimes b)=f(a)b\andeqn \psi(f\otimes b)=f(1)b
$$
for $f\in C(X)$ and $b\in A.$
Since $A$ is an amenable simple \CA\, $\phi$ is a monomorphism.
It follows from \ref{NIPps} that $q_{\omega}(B)$ is a
purely infinite simple \CA. Therefore by \ref{NIL3} there is a sequence of
isometries $s_n\in q_{\omega}(B)$ such that
$$
\lim_{n\to\infty}\|s_n^*\phi(x)s_n-\psi(x)\|=0
$$
for all $x \in C(X)\otimes A.$
It follows from \ref{NIL1} that there is an isometry $s\in q_{\omega}(B)$ such that
$$
s^*\phi(x)s=\psi(x)\,\,\,{\rm for\,\,\, all}\,\,\, x\in C(X)\otimes A.
$$
In particular,
$$
s^*as=s^*\phi(\imath\otimes 1)s=\psi(\imath\otimes 1)=1,
$$
where $\imath$ is the identity function $\imath(t)=t.$
We also have
$$
s^*bs=s^*\phi(1\otimes b)s=\psi(1\otimes b)=b.
$$
Hence $s\in j(A)'$ by Lemma 6.3.6 of \cite{Rob}.
\end{proof}

\begin{thm}\label{NIL5}
Let $B$ be a unital purely infinite simple
\CA\, and $A$ be a unital separable amenable
\CA. Suppose that $\phi, \psi: A\otimes {\cal O}_2\to B$ are two monomorphisms.
Then $\psi$ and $\phi$ are approximately unitarily equivalent.
\end{thm}

\begin{proof}
Let $p_1=\phi(1_{A\otimes {\cal O}_2})$ and $p_2=\psi(1_{A\otimes {\cal O}_2}).$
It follows from 3.6 in \cite{Ro1} that $\phi|_{1\otimes {\cal O}_2}$ and
$\psi|_{1\otimes {\cal O}_2}$ are approximately unitarily equivalent.
It follows that $p_1$ and $q_1$ are equivalent in $B.$ Therefore we may assume, without loss of generality,  that $p_1=p_2.$
 By replacing $B$ by $p_1Bp_1,$ we may further assume that
both $\phi$ and $\psi$ are unital.

Let $\Psi, \Phi: A\otimes {\cal O}_2\to l^{\infty}(B)$ be defined by
$\Psi=\{\psi(a),\psi(a),...,\psi(a),...)$ and
$\Phi(a)=\{\phi(a),\phi(a),...,\phi(a),...)$ for $a\in A,$
respectively.  Fix a free ultrafilter $\omega.$ Put ${\bar
\phi}=\pi\circ \Phi$ and ${\bar \psi}=\pi\circ \Psi,$ where $\pi:
l^{\infty}(B)\to q_{\omega}(B)$ is the quotient map. It follows
from 3.6 in \cite{Ro1} that ${\bar \phi}|_{1\otimes {\cal O}_2}$ and
${\bar \psi}|_{1\otimes {\cal O}_2}$ are approximately unitarily
equivalent. It follows from \ref{NIL1} that they are unitarily
equivalent in $q_{\omega}(B).$ Without loss of generality, we may
assume that
$$
{\bar \phi}|_{1\otimes {\cal O}_2}={\bar \psi}|_{1\otimes {\cal
O}_2}.
$$
Let $D={\bar \phi}(1\otimes {\cal O}_2).$ Then $D\cong {\cal O}_2.$
By  \cite{Ro2} ${\cal O}_2\otimes {\cal O}_2\cong {\cal O}_2.$
Let $\imath: {\cal O}_2\to {\cal O}_2\otimes {\cal O}_2$ be
defined by $\imath(a)=a\otimes 1$ and let $\lambda: {\cal O}_2\otimes {\cal O}_2
\to {\cal O}_2$ be an isomorphism.
Then
$\lambda\circ \imath$ and $\id_{{\cal O}_2}$ are approximately unitarily
equivalent.

Let $\gamma: A\otimes {\cal O}_2\to A\otimes {\cal O}_2\otimes {\cal O}_2
\to A\otimes {\cal O}_2\otimes 1$ be the \hm\, induced by $\lambda\circ \imath$
above. Then ${\bar \phi}\circ \gamma$ is approximately unitarily
equivalent to ${\bar \phi}$ and ${\bar \psi}\circ \gamma$ is
approximately unitarily equivalent to ${\bar \psi}.$
To prove that ${\bar \psi}$ and ${\bar \phi}$ are approximately unitarily
equivalent it suffices to show that ${\bar \phi}\circ \gamma$ and
 ${\bar \psi}\circ \gamma$ are approximately unitarily
equivalent.

There is a unital \SCA\, $D_1 \subset D$ of $q_{\omega}(B)$ which
is isomorphic to ${\cal O}_2$ and its commutant contains  both the
images of ${\bar \phi}\circ \gamma$ and ${\bar \psi}\circ \gamma.$

Let $D_1'$ be the commutant in $q_{\omega}(B).$ Then by \ref{NIP1} $D'$ is
purely infinite simple.
It follows from Corollary \ref{NIL4}
  that there are isometries $s_n, w_n\in D_1'$ such that
$$
\lim_{n\to\infty}\|s_n^*{\bar \psi}\circ \gamma(a)s_n-{\bar \phi}(a)\|=0\andeqn
\lim_{n\to\infty}\|w_n^*{\bar \phi}\circ \gamma(a)w_n-{\bar \psi}(a)\|=0
$$
for all $a\in A.$ Since $((D_1)')'$ contains a unital subalgebra $D_1$ which
is isomorphic to ${\cal O}_2,$ by Lemma 6.3.7 in \cite{Rob},
 ${\bar \phi}\circ \gamma$ and ${\bar \psi}\circ \gamma$ are
approximately unitarily
equivalent. It follows that ${\bar \phi}$ and ${\bar \psi}$ are
approximately unitarily equivalent in $q_{\omega}(B).$
It follows from 6.2.5 in \cite{Rob} that
$\phi$ and $\psi$ are approximately unitarily equivalent.
\end{proof}

\section{Approximately unitarily equivalent extensions}
The purpose of this section is to prove Theorem \ref{IabTTad} and  Theorem \ref{IabT1}. The
statements have been proved for the case that the target \CA\, is a
separable amenable purely infinite simple \CA. The problem we
deal with in this section
is to show a certain absorption property in the absence of
``approximate divisibility" for $M(B)/B.$

\begin{Lem}\label{IabL1}
Let $B$ be a  non-unital and $\sigma$-unital \CA\, and
$A$ be a separable \CA. Let $\tau: A\to M(B)/B$ be
a \hm.
Then there is a sequence of  non-zero mutually orthogonal elements
$a_n\in \tau(A)',$ where $\tau(A)'$ is the commutant of $\tau(A)$
in $M(B)/B.$
\end{Lem}

\begin{proof}
Let $D$ be a separable \CA\, in $M(B)$ such that
$\tau(A)\subset \pi(D),$ where $\pi: M(B)\to M(B)/B$ is the
quotient map. It follows from Lemma 3.1 in \cite{Lnnew} that
there exists an approximate
identity $\{e_n\}$ such that $e_{n+1}e_n=e_ne_{n+1}=e_n,$ $n=1,2,...$ and
$$
\|e_kd-de_k\|\to 0\,\,\,{\rm as}\,\,\, k\to \infty
$$
for all $d\in D.$
Fix a subsequence $X\subset {\mathbb N},$ then $a_X=\sum_{n\in X}(e_{n+1}-e_n)$
is a positive element in $M(A).$ Since $\lim_{k\to\infty}\|e_kd-de_k\|=0$ for each
 $d\in A,$ $\pi(a_X)\pi(d)=\pi(d)\pi(a_X).$ In other words $a_X\in \tau(A)'.$
Suppose that $X$ and $Y$ are two disjoint subset of ${\mathbb N}$ such that
for any $n\in x$ and $m\in Y,$ $|n-m|\ge 2.$ By the assumption that
$e_{n+1}e_n=e_ne_{n+1}=e_n,$ we conclude that $a_Xa_Y=a_Ya_X=0.$
From this it is easy to see that there exists a sequence of nonzero
mutually orthogonal elements in $\tau(A)'.$
\end{proof}

\begin{Lem}\label{IabL2}
Let $A$ be a unital separable amenable \CA, $B$ be a non-unital but
$\sigma$-unital simple \CA\, with a continuous scale  and let
$\omega\in \beta {\mathbb N}\setminus {\mathbb N}$ be a free
ultrafilter. Suppose that $\tau: A\to M(B)/B$ is an essential
unital extension. Let $\tau_{\infty}: A\to l^{\infty}(M(B)/B)$ be
defined by $\tau_{\infty}(a)=(\tau(a),\tau(a),...)$ and let
$\psi=\Phi\circ \tau_{\infty},$ where $\Phi: l^{\infty}(M(B)/B) \to
q_{\omega}(M(B)/B).$ Then there is a unital \SCA\, $C\cong {\cal
O}_{\infty}$ in the commutant of $\psi(A)$ in
$q_{\omega}(M(B)/B).$
\end{Lem}

\begin{proof}
Let $J: M(B)/B\to q_{\omega}(M(B)/B)$ be
defined by $J(b)=\Phi((b,b,...,b,...))$ for $b\in B.$ By
(the proof of ) 7.4 in \cite{Lnpro},
there exists a unital separable purely infinite simple \CA\,
$D\subset M(B)/B$ such that $\tau(A)\subset D.$ It follows from
\ref{IabL1} that there is a sequence of nonzero  mutually
orthogonal positive elements $\{a_n\}$ in $D',$ the commutant of
$D$ in $M(B)/B.$ Let $X=sp(a_1).$ Without loss of generality we
may assume that $\|a_1\|=1$ and $1\in X.$ Define $L_1, L_2:
A\otimes C(X)\to q_{\omega}(M(B)/B)$ by $L_1(x\otimes f)=\psi(x)f(J(a))$
and $L_2(x\otimes f)=\psi(x)f(1)$ for $x\in A$ and $f\in C(X).$ Define
$L_3: D\otimes C(X)\to q_{\omega}(M(B)/B)$ by $L_3(y\otimes
f)=yf(J(a))$ for $y\in D.$ Since $\psi$ is injective and $D$ is
purely infinite simple, one concludes that $L_3$ is injective.
Consequently, $L_1$ is injective. Now we apply an argument in
\cite{KP}. Since $A\otimes C(X)$ is amenable, the Choi-Effros
 lifting theorem
(\cite{CE}) provides unital completely positive lifting $\rho,\sigma:
A\otimes C(X)\to l^{\infty}(B)$ of $L_1$ and $L_2.$
Write
$$
\rho(a)=(\rho_1(a),\rho_2(a),...,\rho_n(a),...) \andeqn
\sigma(a)=(\sigma_1(a), \sigma_2(a),...,\sigma_n(a),...)
$$
for $a\in A\otimes C(X),$ where $\rho_k$ and $\sigma_k$ are
 unital, completely
positive maps from $A\otimes C(X)$ into $B_k.$
It follows from 6.3.5 (iii) in \cite{Rob}
 that there are non-unitary isometries
$s_k\in B_k$ ($k=1,2,...$) such that
$$
\lim_{n\to\infty}\|s_n^*\rho_n(a)s_n-\sigma_n(a)\|=0\,\,\,\,{\rm for\,\,\,
all}\,\,\, a\in A\otimes C(X).
$$
Put $t_1=\pi(s_1,s_2,...,s_n,...))\in l^{\infty}( B).$
Then $t_1$ is a non-unitary isometry.
It follows that
$$
t_1^*\psi(a)t_1=t_1^*L_1(a\otimes 1)t_1=
L_2(a\otimes 1)=\psi(a)
$$
for all $a\in A.$ It follows from  Lemma 6.3.6 in \cite{Rob} that
$t_1\in (\psi(A))'.$ Furthermore,
$$
t_1^*J(a_1)t_1=t_1^*L_1(1\otimes \imath)t_1=L_2(1\otimes \imath)=1,
$$
where $\imath$ is the function $\imath(t)=t.$
Let $t_1t_1^*=q_1.$ Then $q_1\in \psi(A)'$ and $q_1\in
{\overline{J(a_1)\psi(A)'J(a_1)}}.$
We repeat the above argument for $a_2,a_3,....$
Then we obtain a sequence of isometries $t_1,t_2,...$
in $\psi(A)'$
such that $t_n^*t_n=1$ and
$$
\sum_{i=1}^nt_it_i^*\in \overline{(\sum_{i=1}^nJ(a_i) )\psi(A)'
(\sum_{i=1}^nJ(a_i)) }.
$$
It follows that
$\sum_{i=1}^n a_i\le 1.$
Therefore we obtain a unital \SCA\, $C\cong {\cal O}_{\infty}$
in $\psi(A)'.$
\end{proof}

\begin{Lem}\label{Iab3}
Let $A$ be a unital separable amenable \CA\, and $C$ be a unital
separable amenable purely infinite simple \CA. Suppose that
$A\otimes C$ is a unital \SCA\, of a unital \CA\, $B.$ Then there
is an embedding $j: A\otimes C\to A\otimes {\cal O}_2\to B$
satisfying the following:  for any $\ep>0,$ any finite subset
${\cal F}\subset A\otimes C$ and any integer $n>0,$ there exists a
partial isometry $u\in M_{n+1}(B)$ such that $u^*u=1,$
$uu^*=1\oplus j(1_{A\otimes C})\oplus j(1_{A\otimes C})\oplus\cdots
\oplus j(1_{A\otimes C})$ {\rm ( }where $j(1_{A\otimes C})$
repeats $n$ times{\rm )} and
$$
u^*(\id_{A\otimes C}\oplus j\oplus j\oplus \cdots \oplus j)u
\approx_{\ep} \id_{A\otimes C}
\,\,\,\, on\,\,\, {\cal F},
$$
where $j$ repeats $n$ times.
\end{Lem}

\begin{proof}
Let $\ep>0$ and ${\cal G}\subset C$ be a finite subset. It follows
from \cite{KP} that there is \hm\, $\imath: C\to {\cal O}_2\to C$
satisfying the following: for any $\ep>0,$ any finite subset
${\cal G}\subset C$ and any integer $n$ there exists a partial
isometry $w\in M_{n+1}(C)$ such that $w^*w=1_C,$ $ww^*=p=1_C\oplus
\imath(1_C)\oplus \imath(1_C)\oplus \cdots \oplus \imath(1_C)$
(where $\imath (1_C)$ repeats $n$ times) and
$$
w^*(\id_C\oplus \imath\oplus \imath\oplus \cdots \oplus \imath)w\approx_{\ep/2}
\id_C\,\,\,\,\,{\rm on}\,\,\, {\cal G},
$$
where $\imath$ repeats $n$ times.
Let $u=1\otimes w$ and define $j: A\otimes C\to B$ by
$j(a\otimes b)=a\otimes \imath(b))$ for $a\in A$ and $b\in C.$
 One checks that the lemma follows.
\end{proof}

\begin{Lem}\label{Iab5}
Let $A$ be a unital separable amenable \CA\ and $B$ be a non-unital
but $\sigma$-unital simple \CA\, with a continuous scale. Suppose
that $h: A\to M(B)/B$ is a unital injective \hm\, and $\omega\in
\beta{\mathbb N}\setminus {\mathbb N}$ is a free ultrafilter. Let
$\pi:l^{\infty}(M(B)/B)\to q_{\omega}(M(B)/B)$ be the quotient
map. Define $H_0: A\to l^{\infty}(M(B)/B)$ by
$H_0(a)=(h(a),h(a),...,h(a),...)$ and $H=\pi\circ H_0.$ Then,
there exists an injective  \hm\, $j: A\to A\otimes {\cal O}_2
\to  q_{\omega}(M(B)/B)$
satisfying the following: For any $\ep>0,$ any finite subset
${\cal F}\subset A$ and any integer $n>0$ there exist an isometry
$u\in M_{n+1}(q_{\omega}(M(B)/B)$ with $u^*u=1,$ $uu^*=1\oplus
j(1_A)\oplus\cdots \oplus j(1_A)$ {\rm (} where $j(1_A)$ repeats
$n$ times{\rm )} such that
$$
u^*(H\oplus j\oplus j\oplus\cdots \oplus j)u
\approx_{\ep} H\,\,\,\,{\rm on}\,\,\,{\cal F},
$$
where $j$ repeats $n$ times. Moreover, there is $q\in
q_{\omega}(M(B)/B)$ such that $[q]=[H(1_A)]$ and $qj(a)=j(a)q$ for
all $a\in A$ and $qjq$ is an injective full \hm.
\end{Lem}

\begin{proof}
Fix a free ultrafilter $\omega\in \beta{\mathbb N}\setminus
{\mathbb N}.$ We identify $h(A)$ with $A.$ It follows from
\ref{IabL2} that  $H(A)'$ contains a unital \SCA\, $C\cong {\cal
O}_{\infty}.$ Thus we obtain an injective \hm\, $\Psi: A\otimes
{\cal O}_{\infty} \to q_{\omega}(M(B)/B).$ Thus the first part of
the lemma follows from this and \ref{Iab3}.

To prove the very last part of the lemma, we may assume that
$[\imath(1_C)]\not=[1_C],$ where $\imath$ is as in \ref{Iab3}.
There is a projection $q\in C$ such that $q\le \imath(1_C)$ and
$[q]=[1_C].$  Then $qj(a)=j(a)q$ for all $a\in A.$ Since $\Psi$ is
injective, for any nonzero element $a\in A$ and $b\in C,$ $ab=0$
implies that $b=0.$ Thus $qjq$ is injective. To see that $qjq$
is full we note that $q_{\omega}(M(B)/B)$ is purely infinite (see
6.26 in \cite{Rob}).
 \end{proof}

\begin{Def}\label{IIDuct}
{\rm  Let $A$ be a separable amenable \CA\, and $B$ be a
$\sigma$-unital \CA. Then $KL(A,B)$ is defined to be
$KK(A,B)/{\cal T}(A,B),$ where ${\cal T}(A,B)$ is the subgroup  of stable
approximately trivial extensions (see \cite{Lnuct}). When $A$ is
in ${\cal N},$ then $KL(A,B)=KK(A,B)/Pext(K_*(A), K_*(B))$ (see
\cite{Ro3}).

Let $C_n$ be a commutative \CA\, with $K_0(C_n)={\mathbb
Z}/n{\mathbb Z}$ and $K_1(C_n)=0.$ Suppose that $A$ is a \CA. Then
set $K_i(A, {\mathbb  Z}/k{\mathbb  Z})=K_i(A\otimes C_k)$ (see
\cite{Sch1}).  One has the following six-term exact sequence (see
\cite{Sch1}):
$$
\begin{array}{ccccc}
K_0(A) & \to & K_0(A, {\mathbb Z}/k{\mathbb Z}) & \to & K_1(A)\\
\uparrow_{{\bf k}} & & & & \downarrow_{{\bf k}}\\
K_0(A) & \leftarrow & K_1(A,{\mathbb Z}/k{\mathbb Z}) & \leftarrow
&
K_1(A)\,.\\
\end{array}
$$
In \cite{DL}, $K_i(A, {\mathbb Z}/n{\mathbb Z})$ is identified
with $KK^i({\mathbb I}_n, A)$ for $i=0,1$.

 As in \cite{DL}, we use the notation
$$
{\underline K}(A)=\oplus_{i=0,1, n\in {\mathbb  Z}_+}
K_i(A;{\mathbb Z}/n{\mathbb Z}).
$$

By ${ \rm Hom}_{\Lambda}({\underline K}(A),{\underline K}(B))$ we
mean all \hm s from ${\underline K}(A)$ to ${\underline K}(B)$
which respect the direct sum decomposition and the so-called
Bockstein operations (see \cite{DL}). It follows from the
definition in \cite{DL} that if $x\in KK(A,B),$ then the Kasparov
product $KK^{i}({\mathbb I}_n, A)\times x$ gives an element in
$KK^i({\mathbb I}_n, B)$ which we identify with ${ \rm
Hom}(K_i(A,{\mathbb Z}/n{\mathbb Z}), K_0(B,{\mathbb Z}/n{\mathbb
Z})).$ Thus one obtains a map $\Gamma: KK(A,B)\to { \rm
Hom}_{\Lambda}(\underline{K}(A),\underline{K}(B)).$ It is shown by
Dadarlat and Loring (\cite{DL}) that if $A$ is in ${\cal N}$
then, for any $\sigma$-unital \CA\, $B,$ the map $\Gamma$ is
surjective and ${\rm ker}\,\Gamma=Pext(K_*(A),K_*(B)).$ In
particular,
$$
\Gamma: KL(A,B)\to{ \rm Hom}_{\Lambda}({\underline
K}(A),{\underline K}(B))
$$
is an isomorphism. }

\end{Def}

We will use the following theorem.
It is a consequence of the uniqueness theorem
in 5.6.4 in \cite{Lnb} which first appeared in the (preprint)
of \cite{Lnjot}.  It is proved in \cite{Lnuct}.

\begin{thm}{\rm (Theorem 3.9 in \cite{Lnuct})}\label{Q}
Let $A$ be a separable unital amenable \CA\,  and let $B$  a unital
\CA. Suppose that $h_1, h_2: A\to B$ are two unital \hm s such
that
$$
[h_1]=[h_2]\,\,\,\,{\rm in}\,\,\, KL(A,B).
$$
Suppose that $h_0: A\to B$ is  a full unital monomorphism. Then,
for any $\ep>0$ and finite subset ${\cal F}\subset A,$ there is an
integer $n$ and a unitary $W\in U(M_{n+1}(B))$ such that
$$
\|W^*{\rm diag}(h_1(a), h_0(a),\cdots, h_0(a))W -{\rm
diag}(h_2(a), h_0(a),\cdots, h_0(a))\|<\ep
 $$
for all $a\in {\cal F}.$
\end{thm}

\begin{thm}\label{IabTTad}
Let $A$ be a separable amenable \CA\, and $B$ be a non-unital and
$\sigma$-unital simple \CA\, with continuous scale.
Suppose that $\tau_1, \tau_2: A\to M(B)/B$ be two essential extensions.
Then $\tau_1$ and $\tau_2$ are approximately unitarily equivalent
if and only if
$$
[\tau_1]=[\tau_2]\,\,\,{\rm in}\,\,\, KL(A, M(B)/B).
$$
\end{thm}

\begin{proof}
Fix an ultrafilter $\omega\in \beta{\mathbb N}\setminus {\mathbb N}.$
Let $\pi: l^{\infty}(M(B)/B)\to q_{\omega}(M(B)/B)$
denote the quotient
map. Define $\Psi_i: A\to l^{\infty}(M(B)/B)$ by $\Psi_i(A)=
(\tau_i(a),\tau_i(a),...)$
for $a\in A.$ Set $H_i=\pi\circ \Psi_i.$
We will show that for any $\ep>0$ and any finite subset
${\cal F}\subset A$ there is a unitary
$w\in q_{\omega}(M(B)/B)$
such that
$$
{\rm ad} w\circ H_1\approx_{\ep/2} H_2\,\,\,\,{\rm on}\,\,\, {\cal F}.
$$
There are unitaries $u_n\in M(B)/B$
such that $\pi((u_1,u_2,...,u_n,...))=w.$
Therefore
$$
\liminf_{\omega}\|u_n^*h_1(a)u_n-h_2(a)\|\le \ep/2
$$
for $a\in {\cal F}.$
Hence, there exists a subset $X\subset \omega$ such that for all $n\in X,$
$$
\|u_n^*h_1(a)u_n-h_2(a)\|<\ep \,\,\,\,\,{\rm for\,\,\, all}\,\,\, a\in {\cal F}.
$$
Then the theorem follows.

Let $j_1$ and $q$ be as in \ref{Iab5} associated with $H_1$ and $j_2$ be
as in \ref{Iab5} associated with $H_2.$
It follows from  Theorem \ref{Q} that
there is a unitary $z\in
M_{K+1}(q_{\omega}(M(B)/B))$ for some integer
$K>0$ such that
$$
{\rm ad}\,z\circ (H_1\oplus qj_1q\oplus qj_1q\oplus\cdots \oplus qj_1q)
\approx_{\ep/5} H_2\oplus qj_1q\oplus qj_1q\oplus\cdots \oplus qj_1q
\,\,\,{\rm on}\,\,\, {\cal F}.
$$
Therefore (by adding $(1-q)j_1(1-q)\oplus\cdots\oplus (1-q)j_1(1-q)$)
there is a unitary $v\in M_{K+1}(q_{\omega}(M(B)/B))$
such that
$$
{\rm ad}\,v\circ(H_1\oplus j_1\oplus j_1\oplus \cdots \oplus j_1)\approx_{\ep/5}
H_2\oplus j_1\oplus j_1\oplus \cdots \oplus j_1\,\,\,\,{\rm on}\,\,\,{\cal F}.
$$
In particular, we may assume that $v^*(1\oplus j(1_A)\oplus\cdots
\oplus j(1_A))v=1.$
It follows from \ref{Iab5} and \ref{NIL5} that
$$
H_1\oplus j_1\oplus j_1\oplus \cdots \oplus j_1\, {\sim}_{\ep/5}\,
H_1\,\,\,\,{\rm on}\,\,\, {\cal F}.
$$
By \ref{NIL5},
$$
j_1\,{\sim}_{\ep/5}\, j_2.
$$
By \ref{Iab5} and \ref{NIL5} again,
$$
H_2\oplus j_1\oplus j_1\oplus \cdots \oplus j_1\,
{\sim}_{2\ep/5}\, H_2\,\,\,\,{\rm on}\,\,\, {\cal F}.
$$
Combining these inequalities, we obtain
a unitary $w\in q_{\omega}(M(B)/B)$ such that
$$
{\rm ad}\circ w \, H_1\approx_{\ep} H_2\,\,\,{\rm on}\,\,\,{\cal F}.
$$
\end{proof}

When $A$ satisfies the UCT, we have the following approximate
version of Theorem \ref{IabTTad}.
This statement is very close
to that of Theorem 6.3 of \cite{Lnpro}.

\begin{thm}\label{IabTadd}
Let $A$ be a separable amenable \CA\, in ${\cal N}.$
For any $\ep>0$ and finite subset ${\cal F}\subset A,$ there exists
a finite subset ${\cal P}\subset {\underline{K}}(A)$ satisfying the following:
if $h_1, h_2, h_3: A\to C$ are three \hm s, where
$C$ is a unital purely infinite simple \CA\,  such that
$$
[h_1]|_{\cal P}=[h_2]|_{\cal P}
$$
then there is an integer $n>0$ and a unitary $u\in M_{n+1}(C)$
such that
$$
{\rm ad}\circ (h_1\oplus h_3\oplus h_3\oplus \cdots h_3)\approx_{\ep}
h_2\oplus h_3\oplus h_3\oplus \cdots h_3)\,\,\,{\rm on}\,\,\,{\cal F},
$$
where $h_3$ repeats $n$ times.
\end{thm}

\begin{proof}
Let $\{{\cal P}_n\}$ be a sequence of finite subsets of
${\underline{K}}(A)$ such that $\cup_{n=1}^{\infty} {\cal
P}_n={\underline{K}}(A).$ Suppose that there are three sequences
of \hm s $\phi_n, \psi_n, f_n: A\to C_n,$ where $C_n$ is a
sequence of unital purely infinite simple \CA s such that
$$
[\phi_n]|_{{\cal P}_n}=[\psi_n]|_{{\cal P}_n},\,\,\,n=1,2,....
$$
It suffices to show that, there exists $N>0$ and $K>0$ such that
when $n\ge N$ there are unitaries $u_n\in M_{K+1}(C_n)$ satisfying
the following:
$$
{\rm ad} u_n(\phi_n\oplus f_n\oplus f_n\oplus \cdots\oplus
f_n)\approx_{\ep} (\psi_n\oplus f_n\oplus f_n\oplus\cdots\oplus
f_n) \,\,\,\,{\rm on}\,\,\,{\cal F},
$$
where $f_n$ repeats $K$ times.
Let $H_1=\{\phi_n\}, H_2=\{\psi_n\}$ and $H_3=\{f_n\}$ be
\hm s from $A$ into $l^{\infty}(\{C_n\})$ and let
${\bar H_i}=\pi\circ H_i,$ where $\pi: l^{\infty}(\{C_n\})\to
 l^{\infty}(\{C_n\})/c_0(\{C_n\})$ is the quotient map, $i=1,2,3.$
So it suffices to show that there exists $K>0$ such that
there is a unitary $U\in M_{K+1}(l^{\infty}(\{C_n\})/c_0(\{C_n\})$
such that
$$
{\rm ad} U({\bar H}_1\oplus {\bar H}_3\oplus {\bar H}_3\oplus\cdots \oplus
{\bar H}_3)
\approx_{\ep} {\bar H}_2\oplus {\bar H}_3\oplus {\bar H}_3\oplus\cdots \oplus
{\bar H}_3
\,\,\,{\rm on}\,\,\, {\cal F},
$$
where ${\bar H}_3$ repeats $K$ times.

Since each $C_n$ is a unital purely infinite simple \CA, it follows
from 6.5 \cite{Lnuct} that ${\bar H_3}: A\to l^{\infty}(\{C_n\})/c_0(\{C_n\})$
is full. So the theorem follows from
Theorem \ref{Q} if we can show
that
$$
[{\bar H}_1]=[{\bar H}_2]  \,\,\,{\rm in }\,\,\,
KL(A, l^{\infty}(\{C_n\})/c_0(\{C_n\}).
$$
It follows from
Corollary 2.1 in \cite{GL} that, if each $C_n$ is purely infinite
and simple,
$$
K_i(l^{\infty}(\{C_n\}))=\prod_n K_i(C_n),\,\,\,i=0,1,
$$
$$
K_i(l^{\infty}(\{C_n\}),{\mathbb Z}/k{\mathbb Z}))\subset
\prod_n K_i(C_n, {\mathbb Z}/k{\mathbb Z}),i=0,1, k=2,3,...
$$
and
$$
K_i(l^{\infty}(\{C_n\})/c_0(\{C_n\})=\prod_nK_i(C_n)/\oplus_nK_i(C_n), \,\,\,
i=0,1,
$$
$$
K_i(l^{\infty}(\{C_n\})/c_0(\{C_n\}), {\mathbb Z}/k{\mathbb Z}))
\subset \prod_n K_i(C_n, {\mathbb Z}/k{\mathbb Z})/\oplus_n
K_i(C_n, {\mathbb Z}/k{\mathbb Z}),\,\,\, k=2,3,....
$$
Thus, since $[H_1]|_{{\cal P}_n}=[H_2]|_{{\cal P}_n}$ for each $n,$
we conclude from the above computation that
$$
[{\bar H}_1]=[{\bar H}_2] \,\,\,\,{\rm in}\,\,\,
Hom_{\Lambda}({\underline{K}}(A), {\underline{K}}(l^{\infty}(\{C_n\})/c_0(\{C_n\})).
$$
Therefore the theorem follows.
\end{proof}

 The following is an approximate version of Theorem \ref{IabTTad}.

\begin{thm}\label{IabT1}
Let $A$ be a unital separable  amenable \CA\, in ${\cal N}$ and $B$ be a
non-unital but $\sigma$-unital simple \CA\, with a continuous scale.
Suppose that $h_1, h_2: A\to M(B)/B$ are two monomorphisms. For any
$\ep>0$ and any finite subset ${\cal F}\subset A,$ there exists a
finite subset ${\cal P}\subset {\underline{K}(A)}$ satisfying the
following: if
$$
[h_1]|_{{\cal P}}=[h_2]|_{{\cal P}}
$$
then there exists a unitary $u\in M(B)/B$ such that
$$
{\rm ad}\circ h_1\approx_{\ep} h_2\,\,\,\,\,{\it on}\,\,\, {\cal F}.
$$

\end{thm}

\begin{proof}
The proof is exactly the same as that of \ref{IabTTad} but
applying  \ref{IabTadd} instead of \ref{Q}.
\end{proof}

\begin{Cor}\label{Iabcor}
Let $A$ be a separable amenable \CA\, in ${\cal N}$ and $B$ be a
non-unital but $\sigma$-unital simple \CA\, with a continuous scale.
Let $\tau_1, \tau_2: A\to M(B)/B$ be two essential extensions.
Then there exists a sequence of unitaries $u_n\in M(B)/B$ such
that
$$
\lim_{n\to\infty}{\rm ad} u_n\circ \tau_1(a)=\tau_2(a)\,\,\,\,{\it
for\,\,\, all}\,\,\, a\in A
$$
if and only if $[\tau_1]=[\tau_2]$ in $KL(A, M(B)/B).$
\end{Cor}

\section{${\bf Exp}_{ap}(A,B)$}

\begin{Def}\label{IVD1}

{\rm
An essential extension $\tau: A\to M(C)/C$ is said to be {\it
approximately trivial} if there is a sequence of trivial
extensions $\tau_n: A\to M(C)/C$ such that
$\tau(a)=\lim_{n\to\infty}\tau_n(a)$ for all $a\in A.$ Denote by
${\bf Ext}_{ap}(A,B)$ the set of approximately unitarily equivalent classes of
essential extensions. }
\end{Def}

 Let $B$ be  a non-unital but $\sigma$-unital simple \CA\, with a continuous scale
 and $A\in {\cal N}.$
 In this
section we will classify essential extensions of $A$ by $B:$
$$
0\to B\to E\to A\to 0
$$
up to approximately unitary equivalence.

\begin{Lem}\label{IIL2p}
Let $B$ be a unital purely infinite simple \CA\, and let $G_i$ be
a countable subgroup of $K_i(B)$ $(i=0,1).$ There exists a unital
separable purely infinite simple \CA\, $B_0\subset B$ such that
$K_i(B_0)\supset G_i$ and $j_{*i}=\id_{K_i(B_0)},$ where $j:
B_0\to B$ is the embedding.
\end{Lem}

\begin{proof}
Since $B$ is purely infinite, all elements in $K_0(B)$ and in
$K_1(B)$ can be represented by projections and unitaries in $B,$
respectively. Let $p_1,...,p_n,...$ be projections in $B$ and
$u_1,u_2,...,u_n,...$ be unitaries in $B$ such that $\{p_n\}$ and
$\{u_n\}$  generates of $G_0$ and $G_1,$ respectively.
Let $B_1$ be a unital separable purely infinite simple \CA\,
containing $\{p_n\}$ and $\{u_n\}$ (see the proof of 7.4 in
\cite{Lnpro}). Note that
$K_i(B_1)$ is countable. The embedding $j_1: B_1\to B$ gives \hm s
$(j_1)_{*i}: K_0(B_1)\to K_i(B).$ Let $F_{1,i}$ be the subgroup of
$K_0(B_1)$ generated by $\{p_n\}$ and $\{u_n\},$ respectively. It
is clear that $(j_1)_{*i}$ is injective on $F_{1,i},$
$i=0,1.$  In particular, the image of $(j_1)_{*i}$ contains $G_i,$
$i=0,1.$ Let $N_{1,i}'={\rm ker}(j_1)_{*i}$ and
let $N_{1,i}$ be the set of all projections (if $i=0$),
or unitaries (if $i=1$) in $B_1$ which have images in $N_{1, i}'.$
 Let $\{p_{1,n}\}$ be a dense subset of $N_{1,0}$ and
 $\{u_{1,n}\}$ be a dense subset of $N_{1,1},$
respectively. Fix a nonzero projection $e\in
B_1$ such that $[e]=0$ in $K_0(B).$ For each $p_{1,n},$  there
exists a partial isometry $w_{1,n}\in B$ such that
$e=w_{1,n}^*w_{1,n}$ and $w_{1,n}w_{1,n}^*=p_{1,n},$ $n=1,2,....$
For each $u_{1,n},$ there are unitaries $z_{1,n,k}\in B,$ $
k=1,2,...,m(n)$ such that
$$
\|z_{1,n,1}-1\|<1/2, \|z_{1,n,m(n)}-u_{1,n}\|<1/2 \andeqn
\|z_{1,n,k}-z_{1,n,k+1}\|<1/2,
$$
$k=1,2,...,m(n),$ $n=1,2,....$ Let $B_2$ be a separable unital
purely infinite simple \CA\, containing $B_1$ and all
$\{w_{1,n}\}$ and $\{z_{1,n,k}\}.$ Note that if $p\in B_1$ is a
projection and $[p]\in N_{1,0}$ then $[p]=0$ in $K_0(B_2).$
Similarly, if $u\in B_1$ and $[u]\in N_{1,1},$ then $[u]=0$ in
$B_2.$ Suppose that $B_l$ has been  constructed. Let $j_l: B_l\to
B$ be the embedding. Let $N_{l, i}={\rm ker}(j_l)_{*i},$ $i=0,1.$
As before, we obtain a unital separable purely finite simple \CA\,
$B_{l+1}\supset B_l$ such that every projection $p\in B_l$ with
$[p]\in N_{l,0}$ has the property that $[p]=0$ in $K_0(B_{l+1}),$
and every unitary $u\in B_l$ with $[u]\in N_{l,1}$ has the
property that $[u]=0$ in $K_1(B_{l+1}).$ Let $B_0$ be the closure
of $\bigcup_{l=1}^{\infty} B_l.$ Since each $B_l$ is purely infinite
and simple, so is $B_0.$ Note  also that $B_0$ is separable. Let
$j: B_0\to B$ be the embedding.

We claim that
$j_{*i}$ is injective.
Suppose that $p\in B_0$ is a projection such that $[p]\in {\rm ker}j_{*0}$
and $[p]\not=0$ in $B_0.$
Without loss of generality, we may assume that $p\in B_l$ for
some large integer $l.$  Then $[p]$ must be
in the ${\rm ker}(j_l)_{*0}.$  By the construction, $[p]=0$ in
$K_0(B_{l+1}).$ This would imply that $[p]=0$ in $K_0(B_0).$
Thus $j_{*0}$ is injective. An exactly same argument shows that
$j_{*1}$ is also injective.
The lemma then follows.
\end{proof}

\begin{Lem}\label{IILsix}
Let $B$ be a unital purely infinite simple \CA.
Suppose that $G_i\subset K_i(B)$ and
$F_i(k)\subset K_i(B,{\mathbb Z}/k{\mathbb Z})$
are
countable subgroups such that
the image of $F_i(k)\subset K_i(B,{\mathbb Z}/k{\mathbb Z})$ in $K_{i-1}(B)$
is contained in $G_{i-1}$
$(i=0,1, k=2,3,...)$. Then there exists a separable unital
purely infinite simple \CA\, $C\subset B$ such that
$K_i(C)\supset G_i,$ $K_i(C,{\mathbb Z}/k{\mathbb Z})
\supset F_i(k)$ and the embedding $j: C\to B$ induces an
injective map $j_{*i}: K_0(C)\to K_i(B)$ and an injective map
$j_*: K_i(C, {\mathbb Z}/k{\mathbb Z})\to K_i(B,{\mathbb
Z}/k{\mathbb Z}),$
 $k=2,3,....$
\end{Lem}

\begin{proof}
It follows from  \ref{IIL2p} that there is a separable unital
purely infinite simple \CA\, $C_1$ such that $K_0(C_1)\supset G_0$
and $K_1(C_1)\supset G_1$ and $j$ induces an identity map on
$K_0(C_1)$ and $K_1(C_1),$ where $j: C\to B$ is the embedding. Fix
$k,$  let $\{x\in K_i(C_1): kx=0\}=\{g_1^{(i)},g_2^{(i)},...,\}.$
Suppose that $\{s_1^{(i)},s_2^{(i)},...,\}$ is a subset of
$K_i(B,{\mathbb Z}/k{\mathbb Z})$ such that the map from
$K_i(B,{\mathbb Z}/k{\mathbb Z})$ to $K_{i-1}(B)$ maps $s_j^{(i)}$
to $g_j^{(i)}.$ For each $z^{(i)}\in K_i(C_1, {\mathbb
Z}/k{\mathbb Z}),$ there is $s_j^{(i)}$ such that
$z^{(i)}-s_j^{(i)}\in K_i(B)/kK_i(B).$ Since $K_i(C_1)$ is
countable, the set of all possible $z^{(i)}-s_j^{(i)}$ is
countable. Thus one obtains a countable subgroup $G_i^{(')}$ which
contains $K_i(C_1)$ such that $G_i^{(')}/kK_i(B)$ contains the
above the mentioned countable set as well as $F_i(k)\cap (
K_i(B)/kK_i(B))$ for each $k.$ Since countably many countable sets
is still countable, we obtain a countable subgroup
$G_i^{(2)}\subset K_i(B)$ such that $G_i^{(2)}$ contains
$G_i^{(')}$ and $kK_i(B)\cap G_i^{(2)}= kG_i^{(2)},$ $k=1,2,...,$
and $i=0,1.$ Note also $F_i(k)\cap (K_i(B)/kK_i(B))\subset
G_i^{(2)}/kK_i(B).$ By applying \ref{IIL2p}, we obtain a separable
purely infinite simple \CA\, $C_2\supset C_1$ such that
$K_i(C_2)\supset G_i^{(2)}$ and embedding from $C_2$ to $B$ gives
an injective map on $K_i(C_2),$ $i=0,1.$ Repeating what we have
done above, we obtain an increasing sequence of countable
subgroups $G_i^{(n)}\subset K_i(B)$ such that $G_i^{(n)}\cap
kK_i(B)=kG_i^{(n)}$ for all $k$ and $i=0,1$ and an increasing
sequence of separable purely infinite simple \SCA s $C_n$ such
that $K_i(C_n)\supset G_i^{(n)}$ and the embedding from $C_n$ into
$B$ gives an injective map on $K_i(C_n),$ $i=0,1,$ and
$n=1,2,....$ Moreover $F_i^{(k)}\cap (K_i(B)/kK_i(B))\subset
K_i(C_n)/kK_i(B).$ Let $C$ denote the closure of $\bigcup_nC_n$
and $j: C\to B$ be the embedding. Then $C$ is a separable purely
infinite simple \CA\, and $j_{*i}$ is an injective map, $i=0,1.$
We claim that $K_i(C)\cap kK_i(B)=kK_i(C),$ $k=1,2,...,$ and
$i=0,1.$ Note that $K_i(C)=\cup_nG_i^{(n)}.$ Since $G_i^{(n)}\cap
kK_i(B)=kG_i^{(n)} \subset kK_i(C),$ we see that $K_i(C)\cap
kK_i(B)=kK_i(C),$ $i=0,1.$ Thus $K_i(C)/kK_i(C)=K_i(C)/kK_i(B).$
Since $K_i(C)/kK_i(B)\supset F_i^{(k)}\cap(K_i(B)/kK_0(B)),$ we
conclude also that $K_i(C,{\mathbb Z}/k{\mathbb Z})$ contains
$F_i(k).$ Since $j_{*0}$ is injective, $j$ induces an injective
map  from $K_0(C)/kK_0(C)$ into $K_0(B)/kK_0(B)$ for all integer
$k\ge 1.$ Using this fact  and the fact that $j_{*i}: K_i(C)\to
K_i(B)$ is
 injective
  by chasing the following commutative diagram,
 $$
{\small \put(-160,0){$K_0(C)$} \put(0,0){$K_0(C,{\bf Z}/k{\bf
Z})$} \put(180,0){$K_1(C)$} \put(-85,-40){$K_0(B)$}
\put(0,-40){$K_0(B,{\bf Z}/k{\bf Z})$} \put(105,-40){$K_1(B)$}
\put(-85, -70){$K_0(B)$} \put(0,-70){$K_1(B, {\bf Z}/k{\bf Z})$}
\put(105,-70){$K_1(B)$} \put(-160,-110){$K_0(C)$}
\put(0,-110){$K_1(C,{\bf Z}/k{\bf Z})$} \put(180,-110){$K_1(C)$}
\put(-120, 2){\vector(1,0){95}} \put(70,1){\vector(1,0){95}}
\put(-123,-3){\vector(1,-1){30}} \put(30,-3){\vector(0,-1){25}}
\put(180,-2){\vector(-1,-1){30}} \put(-45,-38){\vector(1,0){35}}
\put(70,-38){\vector(1,0){25}} \put(-147, -90){\vector(0,1){85}}
\put(-75,-60){\vector(0,1){15}} \put(115, -45){\vector(0,-1){15}}
\put(190,-7){\vector(0,-1){85}} \put(-7,-68){\vector(-1,0){35}}
\put(95,-68){\vector(-1,0){25}} \put(-123,-102){\vector(1,1){30}}
\put(175, -105){\vector(-1,1){30}} \put(30,-100){\vector(0,1){25
}} \put(-5, -108){\vector(-1,0){100}}
\put(170,-108){\vector(-1,0){95}} \put(-112,-12){$j_{*0}$}
\put(15, -15){$j_* $} \put(150, -14){$j_{*1}$}
\put(-130,-92){$j_{*0}$} \put(15, -88){$j_*$} \put(160,
-88){$j_{*1}$} }
$$
one sees that $j$ induces an injective map from $K_i(C, {\mathbb
Z}/k{\mathbb Z}) $ to $K_i(B, {\mathbb Z}/k{\mathbb Z}).$
\end{proof}

\begin{thm}\label{IIT2}
Let $A$ be a unital separable amenable \CA\, in ${\cal N}$ and $B$
be a non-unital but $\sigma$-unital simple \CA\, with a continuous
scale. Then, for any $x\in KL(A, M(B)/B),$ there exists a
monomorphism  $h: A\to M(B)/B$ such that $[h]=x.$
\end{thm}

\begin{proof}
Put $Q=M(B)/B.$
Since $A$ satisfies the UCT, we may view $x$ as an element in
${\rm Hom}_{\Lambda}(\underline{K}(A), \underline{K}(Q)).$ Note that
$K_i(A)$ is a countable abelian group ($i=0,1$). Let
$G_{0}^{(i)}=\gamma(x)(K_i(A)),$ $i=0,1,$
where $\gamma: {\rm Hom}_{\Lambda}(\underline{K}(A), \underline{K}(Q))
\to {\rm Hom}(K_*(A),K_*(Q))$ is the surjective map. Then $G_{0}^{(i)}$ is
a countable subgroup of $K_i(Q),$ $i=0,1.$
Consider the following commutative diagram:
$$
{\small \put(-160,0){$K_0(A)$} \put(0,0){$K_0(A,{\bf Z}/k{\bf
Z})$} \put(180,0){$K_1(A)$} \put(-85,-40){$K_0(Q)$}
\put(0,-40){$K_0(Q,{\bf Z}/k{\bf Z})$} \put(105,-40){$K_1(Q)$}
\put(-85, -70){$K_0(Q)$} \put(0,-70){$K_1(Q, {\bf Z}/k{\bf Z})$}
\put(105,-70){$K_1(Q)$} \put(-160,-110){$K_0(A)$}
\put(0,-110){$K_1(A,{\bf Z}/k{\bf Z})$} \put(180,-110){$K_1(A)$}.
\put(-120, 2){\vector(1,0){95}} \put(70,1){\vector(1,0){95}}
\put(-123,-3){\vector(1,-1){30}} \put(30,-3){\vector(0,-1){25}}
\put(180,-2){\vector(-1,-1){30}} \put(-45,-38){\vector(1,0){35}}
\put(70,-38){\vector(1,0){25}} \put(-147, -90){\vector(0,1){85}}
\put(-75,-60){\vector(0,1){15}} \put(115, -45){\vector(0,-1){15}}
\put(190,-7){\vector(0,-1){85}} \put(-7,-68){\vector(-1,0){35}}
\put(95,-68){\vector(-1,0){25}} \put(-123,-102){\vector(1,1){30}}
\put(175, -105){\vector(-1,1){30}} \put(30,-100){\vector(0,1){25}}
\put(-5, -108){\vector(-1,0){100}}
\put(170,-108){\vector(-1,0){95}} \put(-111,-14){$\gamma(x)$}
\put(15, -15){$\times x $} \put(148, -14){$\gamma(x)$}
\put(-130,-90){$\gamma(x)$} \put(15, -88){$\times x$} \put(160,
-88){$\gamma(x)$} }
$$
It follows from \ref{IILsix} that there is a unital purely infinite
simple \CA\, $C\subset Q$ such that $K_i(C)\subset G_0^{(i)},$
$K_i(C)\cap kK_i(Q)=kK_i(C),$ $k=1,2,...,$ and $i=0,1,$ and
the embedding $j: C\to Q$ induces injective maps on $K_i(C)$ as well
as on $K_i(C,{\mathbb Z}/k{\mathbb Z})$ for all $k$ and $i=0,1.$
Moreover $K_i(C,{\mathbb Z}/k{\mathbb Z})\supset
(\times x)(K_i(A,{\mathbb Z}/k{\mathbb Z}))$ for $k=1,2,...$ and $i=0,1.$
We
have the following commutative diagram:
$$
{\small \put(-160,0){$K_0(A)$} \put(0,0){$K_0(A,{\bf Z}/k{\bf
Z})$} \put(180,0){$K_1(A)$} \put(-85,-40){$K_0(C)$}
\put(0,-40){$K_0(C,{\bf Z}/k{\bf Z})$} \put(105,-40){$K_1(C)$}
\put(-85, -70){$K_0(C)$} \put(0,-70){$K_1(C, {\bf Z}/k{\bf Z})$}
\put(105,-70){$K_1(C)$} \put(-160,-110){$K_0(A)$}
\put(0,-110){$K_1(A,{\bf Z}/k{\bf Z})$} \put(180,-110){$K_1(A)$}
\put(-120, 2){\vector(1,0){95}} \put(70,1){\vector(1,0){95}}
\put(-123,-3){\vector(1,-1){30}}
\put(180,-2){\vector(-1,-1){30}} \put(-45,-38){\vector(1,0){35}}
\put(70,-38){\vector(1,0){25}} \put(-147, -90){\vector(0,1){85}}
\put(-75,-60){\vector(0,1){15}} \put(115, -45){\vector(0,-1){15}}
\put(190,-7){\vector(0,-1){85}} \put(-7,-68){\vector(-1,0){35}}
\put(95,-68){\vector(-1,0){25}} \put(-123,-102){\vector(1,1){30}}
\put(175, -105){\vector(-1,1){30}}
\put(-5, -108){\vector(-1,0){100}}
\put(170,-108){\vector(-1,0){95}} \put(-112,-12){$\gamma(x)$}
 \put(145, -12){$\gamma(x)$}
\put(-135,-92){$\gamma(x)$}  \put(160, -88){$\gamma(x)$} }
$$
We will add two more maps on the above diagram.
 From
 the fact that  the image of $K_i(A, {\mathbb Z}/k{\mathbb Z})$
 under $\times x$ is contained in $K_i(C,{\mathbb Z}/k{\mathbb
 Z}),$ ($k=2,3,...,$ $i=0,1$),
  we obtain two maps $\beta_i:K_i(A, {\mathbb Z}/k{\mathbb
 Z})\to K_i(C,{\mathbb Z}/k{\mathbb
 Z}),$ $k=2,3,...,$ $i=0,1$ such that $j_*\circ \beta_i=\times x$
 and obtain the
 following commutative diagram:
$$
{\small \put(-160,0){$K_0(A)$} \put(0,0){$K_0(A,{\bf Z}/k{\bf
Z})$} \put(180,0){$K_1(A)$} \put(-85,-40){$K_0(Q)$}
\put(0,-40){$K_0(C,{\bf Z}/k{\bf Z})$} \put(105,-40){$K_1(C)$}
\put(-85, -70){$K_0(C)$} \put(0,-70){$K_1(C, {\bf Z}/k{\bf Z})$}
\put(105,-70){$K_1(C)$} \put(-160,-110){$K_0(A)$}
\put(0,-110){$K_1(A,{\bf Z}/k{\bf Z})$} \put(180,-110){$K_1(A)$}
\put(-120, 2){\vector(1,0){95}} \put(70,1){\vector(1,0){95}}
\put(-123,-3){\vector(1,-1){30}} \put(30,-3){\vector(0,-1){25}}
\put(180,-2){\vector(-1,-1){30}} \put(-45,-38){\vector(1,0){35}}
\put(70,-38){\vector(1,0){25}} \put(-147, -90){\vector(0,1){85}}
\put(-75,-60){\vector(0,1){15}} \put(115, -45){\vector(0,-1){15}}
\put(190,-7){\vector(0,-1){85}} \put(-7,-68){\vector(-1,0){35}}
\put(95,-68){\vector(-1,0){25}} \put(-123,-102){\vector(1,1){30}}
\put(175, -105){\vector(-1,1){30}} \put(30,-100
){\vector(0,1){25}} \put(-5, -108){\vector(-1,0){100}}
\put(170,-108){\vector(-1,0){95}} \put(-112,-12){$\gamma(x)$}
\put(15, -15){$\beta_0 $} \put(145, -12){$\gamma(x)$}
\put(-135,-92){$\gamma(x)$} \put(15, -88){$\beta_1$} \put(160,
-88){$\gamma(x)$} }
$$
Consider the following commutative diagram:
$$
\begin{array}{ccccccc}
\to & K_i(A,{\mathbb Z}/mn{\mathbb Z}) & \to & K_i(A, {\mathbb Z}/n{\mathbb Z})&
\to & K_{i-1}(A, {\mathbb Z}/m{\mathbb Z}) &\to \\
& \downarrow & & \downarrow && \downarrow\\
\to & K_i(Q,{\mathbb Z}/mn{\mathbb Z}) & \to & K_i(Q, {\mathbb Z}/n{\mathbb Z})&
\to & K_{i-1}(Q, {\mathbb Z}/m{\mathbb Z}) &\to \cr
\end{array}
$$
Since $j_*\circ \beta_i=\times x$ and
all vertical maps in the following diagram is injective
$$
\begin{array}{ccccccc}
\to & K_i(C,{\mathbb Z}/mn{\mathbb Z}) & \to & K_i(C, {\mathbb Z}/n{\mathbb Z})&
\to & K_{i-1}(C, {\mathbb Z}/m{\mathbb Z}) &\to \\
& \downarrow & & \downarrow && \downarrow\\
\to & K_i(Q,{\mathbb Z}/mn{\mathbb Z}) & \to & K_i(Q, {\mathbb Z}/n{\mathbb Z})&
\to & K_{i-1}(Q, {\mathbb Z}/m{\mathbb Z}) &\to \cr,
\end{array}
$$
we obtain
the following commutative diagram:
$$
\begin{array}{ccccccc}
\to & K_i(A,{\mathbb Z}/mn{\mathbb Z}) & \to & K_i(A, {\mathbb Z}/n{\mathbb Z})&
\to & K_{i-1}(A, {\mathbb Z}/m{\mathbb Z}) &\to \\
& \downarrow & & \downarrow && \downarrow\\
\to & K_i(C,{\mathbb Z}/mn{\mathbb Z}) & \to & K_i(C, {\mathbb Z}/n{\mathbb Z})&
\to & K_{i-1}(C, {\mathbb Z}/m{\mathbb Z}) &\to \cr
\end{array}
$$
Thus we obtain an element $y\in KL(A, C)$ such that $y\times
[j]=x.$ Since $A$ satisfies the UCT, one checks that $KL(A,
C)=KL(A\otimes {\cal O}_{\infty}, C).$ It follows from 6.6 and 6.7 in
\cite{Lnsemi} that there exists a \hm\, $\phi: A\otimes {\cal
O}_{\infty}\to C\otimes {\cal K}$ such that $[\phi]=y.$
Define
$\psi=\phi|_{A\otimes 1}.$ It is then easy to check that
$[\psi]=y.$ Since $A$ is unital, we may assume that the image
of $\psi$ is in $M_m(C)$ for some integer $m\ge 1.$  Since $C$ is a unital purely infinite simple \CA, $1_m$ is equivalent to a projection in $C.$
Thus we may further assume that $\psi$ maps $A$ into $C.$
 Put $h_1=j\circ \psi.$ To obtain a monomorphism, we
note that there is an embedding $\imath: A\to {\cal O}_2$ (see
Theorem 2.8
in \cite{KP}). Since $M(B)/B$ is purely infinite, we obtain a
monomorphism $\psi: {\cal O}_2\to M(B)/B.$ Let $e=\psi(1_{{\cal
O}_ 2}).$ There is a partial isometry $w\in M_2(M(B)/B)$ such that
$w^*w=1_{M(B)/B}$ and $ww^*=1\oplus e.$ Define $h=w^*(h_1\oplus
\psi\circ \imath)w.$ One checks that $[h]=[h_1]$ and $h$ is a
monomorphism.
\end{proof}

\begin{thm}\label{IITM}
Let $A$ be a separable amenable \CA\, in ${\cal N}$ and $B$ be a
non-unital but $\sigma$-unital simple \CA\, with a continuous scale.
Then there is a bijection:
$$
\Gamma: {\bf Ext}_{ap}(A,B) \to  KL(A, M(B)/B).
$$
\end{thm}

\begin{proof}
This follows immediately from  \ref{Iabcor}.
\end{proof}

\begin{Cor}\label{IIC1}
Let $A$ be a unital separable amenable \CA\, satisfying the UCT and
$B$ be a non-unital but $\sigma$-unital  simple \CA\, with a continuous
scale. Let $\tau$ be a unital essential extension and $\psi: A\to
M(B)$ be a \morp\, such that $\pi\circ \psi=\tau.$ Suppose that
$[\tau]=[t]$ in $KL(A, M(B))$ for some trivial extension $t.$
Then, there exists a sequence of monomorphisms $h_n: A\to M(B)$
such that
$$
\lim_{\to\infty} \pi\circ (h_n(a)-\psi(a))=0
$$
for all $a\in A.$
\end{Cor}

xxxxxxxxxxxxxxxxxxxxxxxxxxxxxxxxxxxxxxx AF xxxxxxxxxxxxxx
\begin{Ex}\label{IIEx}
{\rm Let $B$ be a non-unital separable simple \CA\, with finite
trace and $K_0(B)={\mathbb Q}.$ So $B$ has a continuous scale and
$K_0(M(B)/B)={\mathbb R}/{\mathbb Q}.$  Let $\xi\not=0,1$ in
${\mathbb R}/{\mathbb Q}.$ Suppose that $A$ is a unital separable
amenable \CA\, which satisfies the UCT and suppose that there are
two nonzero elements $g_1, g_2$ in $K_0(A)$ such that $[1_A]=g_1$
and the subgroup generated $g_1$ and $g_2$ is not cyclic. Since
${\mathbb R}/{\mathbb Q}$ is divisible, there is a group \hm\,
$\alpha: K_0(A)\to K_0(M(B)/B)$ such that
$\alpha(g_1)=1$ and $\alpha(g_2)=\xi.$ It follows from
\ref{IITM} that there is an essential unital extension $\tau_{\xi}:
A\to M(B)/B$ such that $(\tau_{\xi})_{*0}=\alpha$ and
$(\tau_{\xi})_{*1}=0.$ Since $K_0(B)={\mathbb Q}$ and $K_1(B)=0,$
we compute that $[\tau_{\xi}]=0$ in $KK(A, B)$ for any such
$\xi.$ However, $[\tau_{\xi}]\not=[\tau_{\xi'}]$ in
$KL(A,M(B)/B)$ if $\xi\not=\xi'.$ This shows that there are
uncountably many non-equivalent essential extensions which
represent the same element in $KK^1(A,B).$
 This example shows how $KL(A, M(B)/B))$ can be used to compute
${\bf Ext}_{ap}(A,B),$ while $KK^1(A,B)$ fails.}
\end{Ex}

\section{Examples}

Theorem \ref{IITM} provides a complete
classification of ${\bf Ext}_{ap}(A,B).$ However, it is not
immediately  clear which elements in $KL(A,M(B)/B)$ give
an approximate trivial extension or a quasidiagonal extension. It
turns out it is rather a complicated problem. First of all from item
(1) below, it could be the case that there are no essential
extensions which are approximately trivial. Second, item (2) and
item (3) below show that
that $[\tau]=0$ in $KL(A,M(B)/B)$ does not imply that $\tau$ is an
approximately
trivial extension.  In this section we will discuss these
problems.

{\it In this section $B$ is a non-unital and $\sigma$-unital
simple \CA\, with real rank zero, stable rank one, weakly unperforated
$K_0(B)$ and with a continuous scale.}

We will show the following:

(1) There are $A$ and $B$ such that there are no trivial essential
extensions of $A$ by $B.$

(2) There are essential extensions $\tau$ such that $[\tau]=0$ in
$KL(A, M(B)/B)$ which are not limits of trivial extensions.

(3) For the same $A$ and $B$ as in (2),  there are
trivial essential extensions $\tau$ such that
$[\tau]\not=0.$


\begin{Ex}\label{IIe1}
{\rm Let $A$ be a unital separable amenable \CA\, and $B$ be a
non-unital but $\sigma$-unital simple \CA\, (with a continuous
scale). It is possible that there are no essential trivial extensions of
the form:
$$
0\to B\to E\to A\to 0.
$$
 For example, let $A={\cal O}_n$  ($n\ge 2$) and $B$
be any non-unital AF-algebra with a continuous scale. There are many
extensions of $A$ by $B.$ This is because $M(B)/B$ is purely
infinite simple and one can easily find monomorphisms from ${\cal
O}_n$ into $M(B)/B.$ But none of them are splitting. In fact there
is no monomorphism $h: A\to M(B).$ Since $M(B)$ admits a tracial
state, $h(A)$ would have a tracial state too. But this is
impossible.

   From this example, one sees clearly that for many
\CA s $A$ there is no single  essential trivial extension of $A$
by $B.$ Therefore some restriction on $A$ is
needed to guarantee that there are trivial essential extensions.
}

\end{Ex}

\begin{Lem}\label{IIILemb}
Let $A$ be a unital  AF-algebra such that there is a positive
\hm\, $\alpha: K_0(A)\to Aff(T).$ Then there exists a \hm\, $h:
A\to M(B)$ such that $h_{*0}=\alpha$ and $h(A)\cap B=\{0\}.$
\end{Lem}

\begin{proof}
It is easy to see and known that the lemma holds for the case
that $A$ is finite
dimensional.
Let $F$ be a finite dimensional \CA.
Suppose that $h_1, h_2: F\to M(B)$ are two \hm s such that
$h_i(F)\cap B=\{0\}.$ Suppose also that $(h_1)_{*0}=(h_2)_{*0}.$
Then by (2) in \ref{0T1},
it is standard to see that
$h_1$ and $h_2$ are unitarily equivalent.

Now let $A$ be the closure of $\bigcup_{n=1}^{\infty} A_n,$ where
$A_n\subset A_{n+1}$ and ${\rm dim} A_n<\infty.$
Denote by $j_n$ the embedding from $A_n$ to $A.$
Let $\alpha_n=\alpha\circ (j_n)_{*0}.$
Let $h_1: A_1\to M(A)$ be such that $h_1(A_1)\cap B=\{0\}$ and
$(h_1)_{*0}=\alpha_1.$
Suppose that $h_m: A_m\to M(B)$ has been defined
such that $h_m|_{A_j}=h_j$ for $j<m,$ $h_m(A_m)\cap B=\{0\}$ and
$(h_m)_{*0}=\alpha_m.$
Let $\phi_{m+1}: A_{m+1}\to M(B)$ be such that
$\phi_{m+1}(A_{m+1})\cap M(B)=\{0\}$ and $(\phi_{m+1})_{*0}
=\alpha_{m+1}.$
Let $\imath_{m}: A_m\to A_{m+1}$ be the embedding.
Then we have
$\alpha_m=(h_m)_{*0}=(\phi_{m+1}\circ \imath_{m}).$
From what we have shown there is a unitary $u_{m+1}\in M(B)$ such
that
$$
{\rm ad}\circ \phi_{m+1}\circ \imath_{m}=h_m.
$$
Put $h_{m+1}={\rm ad}\circ \phi_{m+1}.$ Then we have the following
commutative diagram:
$$
\begin{array}{ccccccccc}
A_1 &\to &A_2& \to & A_3& \to &\cdots & A\\
\downarrow_{h_1} & &\downarrow_{h_2} &&\downarrow_{h_3} &\cdots\\
M(B)&{\stackrel{\id_M(B)}{\to}} & M(B)&{\stackrel{\id_M(B)}{\to}}&
M(B) & {\stackrel{\id_M(B)}{\to}} &\cdots &M(B)\cr
\end{array}
$$
It follows that there is a monomorphism $h: A\to M(B)$ such
that $h(A)\cap B=\{0\}.$
\end{proof}

\begin{thm}\label{IIITtriv}
Let $A$ be a separable amenable \CA\, satisfying the UCT. Suppose
that $A$ can be embedded into a unital simple AF-algebra. Then for
any $B$ there exists an essential trivial extension $\tau$ of $A$
by $B.$
\end{thm}

\begin{proof}
Suppose that $C$ is a unital simple AF-algebra and $j: A\to C$ is
an embedding. Let $t$ be a normalized trace on $C.$ Define $\beta:
K_0(C)\to Aff(T)$ by  $\beta([p])=t(p)[1_{M(B)/B}]$ for projection
$p\in C.$ Then $\beta$ is a positive \hm. It follows from
\ref{IIILemb} that there is a monomorphism $h: C\to M(B)$ such
that $h_{*0}=\beta$ and $h(C)\cap B=\{0\}.$  Define $\phi: A\to
M(B)$ by $\phi=h\circ j.$ One sees that $\phi$ give an essential
trivial extension of $A$ by $B.$
\end{proof}

Suppose that there are trivial essential extensions
of $A$ by $B.$
One  would like to know
when an extension trivial, or when
an extension is the limit of trivial extensions.

\begin{Ex}\label{IIe2}
There are essential extensions $\tau$ which are not approximately
trivial but  $[\tau]=0$ in\\
$KL(A,M(B)/B).$

 {\rm Let $A$ be a unital separable amenable \CA\, and let
$\tau$ be an essential extension of $A$ by $B$ such that
$[\tau]=0$ in $KL(A,M(B)/B).$ Such $\tau$ exists (by \ref{IITM} or by
first mapping $A$ to ${\cal O}_2$ and then mapping ${\cal O}_2$
into $M(B)/B$).

To be more precise, we let $A$ be the unital simple AF-algebra
with $K_0(A)=D_{\theta},$ where $\theta$ is  an irrational
number and
$$
D_{\theta}=\{ m+n\theta: m, n\in {\mathbb Z}\}.
 $$
 with usual
order inherited  from ${\mathbb R}. $ We may assume that
$[1_A]=1.$ Let $B$ be a non-unital (non-zero) hereditary \SCA\, of
the UHF-algebra with $K_0(B)={\mathbb Z}[1/2].$ Note that $B$ has a
unique normalized trace. So it has a continuous scale. We further
assume that $[1_{M(B)/B}]=0.$ So there is an essential extension
$\tau$ of $A$ by $B$ such that $[\tau]=0.$
 However, there is no (non-zero) positive \hm\,
$\alpha: K_0(A)\to K_0(B).$ If there is a  trivial extension
$\tau$ of $A$ by $B$ with $[\tau]=0$ in $KL(A,M(B)/B),$ then
$\tau_{*0}: K_0(A)\to K_0(M(B)/B)$ is zero. It follows \ref{IIILemb} that
$$
K_0(M(B)/B)=Aff(T(B))/K_0(B)={\mathbb R}/{\mathbb Z}[1/2].
$$
If $\tau$ is trivial, there would be a monomorphism $h: A\to M(B)$
such that $h|_{*0}$ maps $K_0(A)$ to $K_0(B)\subset Aff(T(B))$
positively. However there is no positive \hm\, from $D_{\theta}$
into ${\mathbb Z}[1/2].$ In fact any positive \hm\, from
$D_{\theta}$ into ${\mathbb R}$ has to be the form (see \ref{0Pph})
$$
h_{*0}(r)=(h)_{*0}(1)r \,\,\,{\rm for\,\,\, all }\,\,\, r\in
D_{\theta}.
$$
 So $\tau$ can never be trivial.
Furthermore, $\tau$ cannot be approximately trivial. To see this,
assume that $\tau_n: A\to M(B)/B$ are trivial extensions such that
$$ \lim_{n\to\infty}\tau_n(a)=\tau(a)$$
for  all $a\in A.$  Let $G_0\subset K_0(A)$ which contains $1$ and
$\theta.$ Thus, for all large $n,$ $(\tau_n)_{*0}(\theta)=0.$
Suppose that $h_n: A\to M(B)$ such that $\pi\circ h_n=\tau_n,$
$n=1,2,...,$ where $\pi: M(B)\to M(B)/B$ is the quotient map. Thus
$(h_n)_{*0}$ is a positive \hm\, from $D_{\theta}$ into ${\mathbb
R}.$ From the above expression of $h_{*0}$ (see \ref{0Pph}) we see that
$(h_n)_{*0}$   cannot map both $1$ and $\theta$ into rational
numbers. In other words, such an $h_n$ does not exist. Hence $\tau$
is not approximately trivial.}
\end{Ex}

\begin{Ex}\label{IIIEx0not}
Nevertheless, there are essential trivial extensions of $A$ by
$B$ such that $[\tau]\not=0$ in $KL(A,M(B)/B).$

{\rm  Let $s$ be the unique normalized trace on $B.$ Suppose that
$[1_A]=1$ in $D_{\theta}.$ Let $\beta: D_{\theta}\to {\mathbb R}$
the usual embedding. It follows from \ref{IIILemb} that there is a
monomorphism $h: A\to M(B)$ such that $h_{*0}=\beta$ and $h(A)\cap
B=\{0\}.$  Let $\tau=\pi\circ h,$ where $\pi: M(B)\to M(B)/B$ be
the quotient map. Then $\tau$ is a trivial essential extension.
However, $\tau_{*0}: {\mathbb Q}\to {\mathbb R}/{\mathbb Z}[1/2]$
($\, \cong K_0(M(B)/B))$ is not zero. Therefore $[\tau]\not= 0$ in
$KL(A, M(B)/B).$}

\end{Ex}

\section{Quasidiagonal extensions --- general and infinite cases}

\begin{Def}\label{IVDquas}
{\rm  Let $A$ be a separable \CA, $C$ be a non-unital but
$\sigma$-unital \CA\, and $\tau: A\to M(C)/C$ be  an essential
extension. Let $\pi: M(C)\to M(C)/C$ be the quotient map. Set
$E=\pi^{-1}(\tau(A)).$ The extension $\tau$ is  said to be {\it
quasidiagonal} if there exists  an approximate identity $\{e_n\}$
of $C$ consisting of projections such that
$$
\lim_{n\to\infty}\|e_nb-be_n\|=0
$$
for all $b\in E.$

Suppose that there is a bounded linear map
 $L: A\to M(B)$
 such that $\pi\circ L=\tau.$ Then
$$
\|e_nL(a)-L(a)e_n\|\to 0\,\,\,\, {\rm as}\,\,\,n\to\infty
$$
for all $a\in A.$ }

\end{Def}

In this section and the next, we will study quasidiagonal extensions.
The first question is when  quasidiagonal extensions exist.

\begin{thm}\label{IVqkadd}
Let $A$ be a separable amenable \CA\, and $B$ be a non-unital and
$\sigma$-unital \CA.
Suppose that $\tau: A\to M(B)/B$ is an essential quasidiagonal extension.
Then for each finitely generated subgroup $G$ of $\underline{K}(A)$
there exists a \hm\, $\alpha: G\to \underline{K}(M(B))$ such that
$$
\pi_*\circ \alpha|_G=(\tau_*)|_{G},
$$
where $\pi: M(B)\to M(B)/B$ is the quotient map.
\end{thm}

\begin{proof}
Suppose that $\tau: A\to M(B)/B$ is a quasidiagonal essential
extension
 of $A$ by $B.$ Let $L: A\to M(B)$ be a \morp\,
 such that $\pi\circ L=\tau.$
Since $\tau$ is quasidiagonal, there exists an approximate
identity $\{e_n\}$ for $B$ such that
$$
\|e_nL(a)-L(a)e_n\|\to 0,\,\,\,{\rm as}\,\,\, n\to\infty
$$
for all $a\in A.$
Define $L_n: A\to M(A)$ by $L_n(a)=(1-e_n)L(a)(1-e_n)$ for
$a\in A.$
Then $\{L_n\}$ is a sequence of asymptotically
multiplicative \morp s.
It follows that for each finitely generated subgroup $G$ of
$\underline{K}(A),$ there exists $N>0$ such that
 $\{L_n\}$ gives a \hm\, $\alpha_n: G\to
\underline{K}(M(B))$ for all $n\ge N.$
 Since $\pi\circ L_n=\tau$ for all $n,$
it follows that $\pi_*\circ \alpha|_{G}=(\tau_*)|_{G}.$
\end{proof}

\begin{Cor}\label{IVCadd1}
Let $A$ be a separable amenable \CA\, and $B$ be a $\sigma$-unital
and stable \CA.
Suppose that $\tau: A\to M(B)/B$ is an essential quasidiagonal extension.
Then $\tau$ induces the zero map from $\underline{K}(A)$
to $\underline{K}(M(B)/B).$
Furthermore the six-term exact sequence in $K$-theory associated with
the extension splits into two pure extensions of groups:
$$
0\to K_i(B)\to K_i(E)\to K_i(A)\to 0\,\,\,\, i=0,1.
$$
\end{Cor}

\begin{proof}
This follows from the fact that when $B$ is stable,
$K_i(M(B))=0$ ($i=0,1$) and \ref{IVqkadd}.
\end{proof}

\begin{Lem}\label{IVLqalim}
Let $A$ be a separable amenable \CA\, and $C$ be a
non-unital but $\sigma$-unital \CA.
Suppose that $\tau: A\to M(C)/C$ is an essential extension
such that there exists a sequence of quasidiagonal extensions
$\tau_n: A\to M(C)/C$ such that
$$
\lim_{n\to\infty} \tau_n(a)=\tau(a)
$$
for all $a\in A.$ Then $\tau$ is quasidiagonal.
\end{Lem}

\begin{proof}

Let $\{a_n\}$ be a dense sequence in  the unit ball of $A.$
Suppose that
$$
\|\tau_n(a)-\tau(a)\|<1/2^{n+3}
$$
for all $a\in \{a_1,a_2,...,a_n\},$ $n=1,2,....$
Let $L_n: A\to M(B)$ be a \morp\,
such that $\pi\circ L_n=\tau_n.$
There exists an approximate identity $\{e_k^{(n)}\}$ for $B$
consisting of projections such
that
$$
\lim_{n\to\infty}\|e_k^{(n)}L_n(a)-L_n(a)e_k^{(n)}\|=0
$$
for all $a\in A.$
Let $L: A\to M(B)$ be a \morp\, such that
$\pi\circ L=\tau.$
Suppose that $b$ is a strictly positive element
for $B.$
We may assume
that
$$
\|(1-e_1^{(1)})(L_1(a)-L(a))\|<1/2^{3}\andeqn
\|(L_1(a)-L(a))(1-e_1^{(1)})\|<1/2^{3}
$$
for $a=a_1.$ Put $q_1=e_1^{(1)}.$
Note that
$$
\|q_1L(a_1)-L(a_1)q_1\|<1/2.
$$
By changing notation if necessary, we may assume that
$$
\|(1-e_2^{(2)})(L_2(a)-L(a))\|<1/2^{2+2}\andeqn
\|(L_2(a)-L(a))(1-e_2^{(2)})\|<1/2^{2+2}
$$
for $ a\in \{a_1,a_2\}$ as well as
$$
\|(1-e_2^{(2)})q_1\|<1/2^{2+3}\andeqn \|(1-e_2^{(2)})b\|<1/2^3.
$$
There is a projection $q_2\ge q_1$ such that
$$
\|e_2^{(2)}-q_2\|<1/2^3.
$$
Note also
$$
\|q_2L(a)-L(a)q_2\|<1/2^2\,\,\,{\rm for}\,\,\, a\in \{a_1,a_2\}.
$$
We also have
$$
\|(1-q_2)b\|<1/2^2.
$$
We may assume that
$$
\|(1-e_3^{(3)})(L_3(a)-L(a))\|<1/2^{3+2} \andeqn
\|(L_3(a)-L(a))(1-e_3^{(3)}\|<1/2^{3+2}
$$
for $a\in \{a_1,a_2,a_3\}$ as well as
$$
\|(1-e_3^{(3)})q_i\|<1/2^{3+3}, i=1,2 \andeqn
\|(1-e_3^{(3)})b\|<1/2^4.
$$
There exists a projection $q_3\ge q_2$ such that
$$
\|e_3^{(3)}-q_3\|<1/2^4.
$$
Thus  we have
$$
\|q_3L(a)-L(a)q_3\|<1/2^3\,\,\,{\rm for}\,\,\, a\in \{a_1,a_2,a_3\}.
$$
We also have
$$
\|(1-q_3)b\|<1/2^3.
$$
We continue in this fashion. It follows that we obtain an
increasing sequence of projections $\{q_n\}$ in $B$ such that
$$
\|q_nL(a)-L(a)q_n\|<1/2^n\,\,\,\, a\in \{a_1,a_2,...,a_n\}
$$
and
$$
\|(1-q_n)b\|<1/2^n,
$$
$n=1,2,....$
It remains to show that $\{q_n\}$ is an approximate identity
for $B.$
Since $\lim_{n\to\infty}\|(1-q_n)b\|=0,$ one concludes
that for any positive function $f\in C_0((0,\|b\|]),$
$$
\lim_{n\to\infty}\|(1-q_n)f(b)\|=0.
$$
For any $a\in A$ and $\ep>0$ there exists a positive function
$f\in C_0((0,\|b\|])$ such that
$$
\|f(b)a-a\|<\ep/3.
$$
Choose $N>0$ such that
$$
\|(1-q_n)f(b)\|<\ep/3(\|a\|+1)\,\,\,{\rm for \,\,\, all}\,\,\, n\ge N.
$$
Then, for $n\ge N,$
$$
\|(1-q_n)a\|\le \|q_na-q_nf(b)a\|+\|q_nf(b)a-f(b)a\|+\|f(b)a-a\|
<\ep.
$$
It follows that $\{q_n\}$ is an approximate identity for $B.$
\end{proof}
Results in this paper can be also used to prove the following.

\begin{thm}{\rm (Brown-Salinas-Schochet)}\label{IVTbss}
Let $A$ be a separable amenable \CA\, in ${\cal N}$ and $B$ be a $\sigma$-unital
stable \CA. Suppose that there exists an essential quasidiagonal extension of
$A$ by $B.$ The zero element in $KL(A, M(B)/B)$ corresponds to
the set of stably quasidiagonal extensions as well as
stably approximately trivial extensions.
\end{thm}

\begin{proof}
It is proved in \cite{Lnuct} that stably approximately trivial extensions
correspond to the zero element in $KL(A, M(B)/B)$ without assuming $A$ satisfies
the UCT. Corollary \ref{IVCadd1} proves that quasidiagonal extensions give a zero
element in $KL(A, M(B)/B).$ It follows from 3.9 in \cite{Lnuct} that
extensions which represent the same element in $KL(A,M(B)/B)$ are stably
approximately unitarily equivalent. Then by \ref{IVLqalim}
every (stably) approximately trivial extension is stably quasidiagonal.
\end{proof}

\begin{thm}\label{IVTinq}
Let $A$ be a separable exact \CA\, and $B$ be a $\sigma$-unital
purely infinite simple \CA.
Then there are essential quasidiagonal extensions
\end{thm}

\begin{proof}
Let $e\in B$ be a nonzero projection such that $[e]=0$ in $K_0(B).$
Then by \cite{Br1} $eBe\otimes {\cal K}\cong B\otimes {\cal K}.$
It follows from a result of S. Zhang that $B\cong B\otimes {\cal K}$
(\cite{Zh3}).
Thus we obtain an approximate identity $\{e_n\}$ of $B$ such that
each $e_n$ is a projection and $[e_n]=0.$
Since $A$ is exact, by Theorem 2.8 in \cite{KP}, there
exists a monomorphism $\imath: A\to {\cal O}_2.$
Since $e_nBe_n$ is purely infinite and $[e_n]=0,$ there
is an embedding $\phi_n: {\cal O}_2\to e_nBe_n.$
Now define
$$
\psi(a)=\sum_{n=1}^{\infty} \phi_n\circ \imath(a)\,\,\,{\rm for}\,\,\,
a\in A.
$$
Then $\psi$ is an injective \hm\, from $A$ into $M(B)$ such that
$\psi(A)\cap B=\{0\}.$
Let $\tau=\pi\circ \psi.$ Then $\tau$ is an essential quasidiagonal
extension of $A$ by $B.$
\end{proof}

The following
follows from \ref{IVTinq}, Theorem 1.4 in \cite{Sch2} and Kirchberg's
absorbing theorem \cite{K1}.
It also follows from \ref{IVTinq} and  \ref{IVCadd1}.

\begin{Cor}\label{IVTqsin}
Let $A$ be a separable amenable \CA\, in ${\cal N}$
 and $B$ be a non-unital and $\sigma$-unital
purely finite simple \CA.
Suppose that $\tau: A\to M(B)/B$ is an essential  extension.
Then $\tau$ is quasidiagonal extension if and only if it
is an approximately trivial extension, and,
 if and only if
$\tau$ induces a zero element in $KL(A, M(B)/B).$
\end{Cor}

In the next section we will discuss the case that $B$ is not purely infinite.

\section{Quasidiagonal extensions -- finite case}

\begin{Def}\label{IVDcq}
{\rm Recall that a separable \CA\, is said to be quasidiagonal if
there exists a faithful representation $\phi: A\to B(H)$ for some
separable Hilbert space $H$ such that
$$
\|p_n\phi(a)-\phi(a)p_n\|\to  0\,\,\,{\rm as}\,\,\,n\to\infty
$$
for all $a\in A,$ where $\{p_n\}$ is an approximate identity of
${\cal K}$ consisting of finite rank projections.

All AF-algebras are quasidiagonal. All commutative \CA s are
quasidiagonal. All AH-algebras are quasidiagonal.
All residually finite dimensional \CA s are quasidiagonal.
Inductive limits of quasidiagonal \CA s are quasidiagonal.
}

\end{Def}

Recall that a \CA\, $A$ is said to have the property (SP)
if every non-zero hereditary \SCA\, contains a nonzero projection.
One should note that every \CA\, with real
rank zero has the property (SP) but the converse is not true.

\begin{thm}\label{IVTqa}
Let $A$ be a separable quasidiagonal amenable \CA\, and $C$
be a non-unital but $\sigma$-unital simple \CA\, which
admits an approximate identity consisting of projections and
has the
property (SP).
Then there
exists an (essential) quasidiagonal extension $\tau: A\to M(C)/C.$
\end{thm}

\begin{proof}
We may assume that $C\not={\cal K}.$
There is a sequence of \morp s $L_n: A\to F_n,$ where $F_n$ are
finite dimensional \CA s, such that
$$
\|L_n(ab)-L_n(a)L_n(b)\|\to 0,\,\,\,{\rm as}\,\,\, n\to\infty
$$
for all $a, b\in A.$
Let $\{a_n\}$ be a dense sequence in the unit ball of $A.$
By passing to a subsequence, we may assume that
$$
\|L_n(ab)-L_n(a)L_n(b)\|<1/2^{n+1}
$$
for all $a, b\in \{a_1,...,a_n\}.$

Let $\{e_n\}$ be an approximate identity for $B$ consisting of projections.
We may assume that $e_{n+1}-e_n\not=0$ for all $n.$
It is known (see 3.5.7 in \cite{Lnb}) that there is a monomorphism
$j_n: F_n\to (e_{n+1}-e_n)C(e_{n+1}-e_n).$
Define
$$
L(a)=\sum_{n=1}^{\infty} j_n\circ L_n(a)\,\,\,{\rm for}\,\,\, a\in A.
$$
Note that the sum converges in the strict topology. One checks
that $L: A\to M(C)$ is a (completely) positive linear
contraction. Note also, for any $a, b\in \{a_1,a_2,...,a_n\},$
$$
\|\sum_{k=n}^{n+m} (L_n(ab)-L_n(a)L_n(b))\|<1/2^{n}\,\,\,{\rm for\,\,\, all}
\,\,\, m>0.
$$
This implies
$$
L(ab)-L(a)L(b)=\sum_{n=1}^{\infty} j_n\circ (L_n(ab)-L_n(a)L_n(b))\in C
$$
for all $a, b\in A.$
Let $\tau=\pi\circ L,$ where $\pi: M(C)\to M(C)/C$ is the quotient
map. Then $\tau$ is an essential quasidiagonal extension.

\end{proof}

\begin{thm}\label{IVTqa2}
Let $A$ be a separable amenable \CA\, and let  $C$ be as in
\ref{IVTqa}.  Suppose that, in addition, $C$ is also a
quasidiagonal \CA. Then there is an (essential) quasidiagonal
extension of $A$ by $B$ if and only if $A$ is quasidiagonal.
\end{thm}

\begin{proof}
It suffices to show the ``only if" part.
Suppose that $L: A\to M(C)$ is a bounded linear map
such that $\pi\circ L: A\to M(C)/C$ is a monomorphism
and
$$
\lim_{n\to\infty}\|e_nL(a)-L(a)e_n\|=0\,\,\,{\rm for\,\,\, all}\,\,\, a
\in A,
$$
where $\{e_n\}$ is an approximate identity consisting
of projections.
Let $\{a_n\}$ be a  dense sequence of $A.$ We may assume  that
$$
\|e_nL(a)-L(a)e_n\|<1/2^n
\,\,\,{\rm for}\,\,\,a\in \{a_1,a_2,...,a_n\},
$$
$n=1,2,....$
Since $\pi\circ L$ is a monomorphism,
we may further assume that
$$
\|e_nL(a)e_n\|\ge \|a\|-1/2^n\,\,\,{\rm for}\,\,\, a\in \{a_1,a_2,...,
a_n\},
$$
$n=1,2,....$
Since $C$ is quasidiagonal, it follows that
there exists a finite dimensional \CA\, $F_n$ and
a \morp\, $\phi_n: e_nCe_n\to F_n$ such that
$$
\|\phi_n(b)\|\ge \|b\|-1/2^n\andeqn
\|\phi_n(bc)-\phi_n(b)\phi_n(c)\|<1/2^n
$$
for $b\in e_nL(a_i)e_n,$ $i=1,2,...,n$ and $n=1,2,....$\,
Define $L_n: A\to F_n$ by $L_n(a)=\phi_n(e_nL(a)e_n)$ for $a\in A.$
Then
$$
\|L_n(a_i)\|\ge \|a_i\|-1/2^{n-1}, i=1,2,...,n \andeqn
\lim_{n\to\infty}\|L_n(ab)-L_n(a)L_n(b)\|=0
$$
for all $a, b\in A.$ It follows from Theorem 1 in \cite{V2} that
$A$ is quasidiagonal.
\end{proof}

\vspace{0.2in} {\it For the rest of this section, $B$ is always a
non-unital but $\sigma$-unital simple \CA\, with real rank zero,
stable rank one, weakly unperforated $K_0(B)$ and a continuous
scale.}

\vspace{0.2in}

\begin{Lem}\label{II1Llift}
Let $A$ be a finite dimensional \CA.  Let $\tau: A\to M(B)/B$ be a
monomorphism such that ${\rm im}\tau_{*0}\subset
Aff(T)/\rho_B(K_0(B)).$ Then $\tau$ is trivial and there is a
monomorphism $h: A\to M(B)$ such that $\pi\circ h=\tau.$
\end{Lem}

\begin{proof}
Suppose that $A=M_{r(1)}\oplus\cdots M_{r(k))}.$ So $K_0(A)$ is
$k$ copies of ${\mathbb Z}.$ Let $e_i$ be a minimal projection in
$M_{r(i)},$ $i=1,2,...,k.$ There are $x_i\in
Aff(T(B))/\rho(K_0(B))$ such that $[\tau(e_{i})]=x_i,$
$i=1,2,...,k$ and $\sum_{i=1}^{k}r(i)x_i=[\tau(1_A)].$ It
follows from Lemma 1.3 in \cite{Lnamj} that there are projections
$q_i\in M(B)$ such that $[\pi(q_i)]=r(i)[\tau(e_i)].$ Thus we
obtain a positive \hm\, $\alpha: K_0(A)\to Aff(T)$ such that
$\pi_{*0}\circ \alpha=\tau_{*0}.$ It follows from \ref{IIILemb} that
there is a monomorphism $h: A\to M(B)$ such that $h(A)\cap
B=\{0\}$ and $h_{*0}=\alpha.$ It then follows that $ [\pi\circ
h]=[\tau]$ in $KL(A, M(B)/B).$ Since $A$ is finite dimensional, it
follows that $\pi\circ h$ is  unitarily equivalent to
$\tau.$ This implies that $\tau$ is trivial.
\end{proof}

\begin{Lem}\label{IVLqanes}
Let $A$ be a separable amenable  \CA. Suppose that
$\tau: A\to M(B)/B$ is an essential quasidiagonal extension. Then
$\tau_{*1}=0,$  ${\rm im}\tau_{*0}\subset Aff(T)/\rho_B(K_0(B))$
and $[\tau]|_{K_i(A, {\mathbb Z}/k{\mathbb Z})}=0$ for $i=0,1$ and
for all $k\ge 2.$
\end{Lem}

\begin{proof}
Since $K_1(M(B))=\{0\}$ (see \ref{0T1}),  by \ref{IVqkadd}, $\tau_{*1}=0.$
By \ref{IVqkadd} and \ref{0T2}, ${\rm im}\tau_{*0}\subset Aff(T)/\rho_B(K_0(B)).$
Since $K_0(M(B))=Aff(T)$ is torsion free, $\tau_{*0}|_{{\rm tor}(K_0(A))}=0.$
Let $C_k$ be as in \ref{IIDuct}. Then  $L\otimes \id_{C(C_k)}: A\otimes
C(C_k)\to M(B)\otimes C(C_k)$ lifts $\tau\otimes \id_{C(C_k)}.$
Since $K_0(M(B))=Aff(T)$ is divisible and $K_1(M(B))=0,$
$K_0(M(B),{\mathbb Z}/k{\mathbb Z})=\{0\}$ for all $k\ge 2.$
It follows from \ref{IVqkadd} that
$[\tau]|_{K_0(A, {\mathbb Z}/k{\mathbb Z})}=0.$

Note also that since $K_0(M(B))$ is torsion free and
$K_1(M(B))=0,$ $K_1(M(B),{\mathbb Z}/k{\mathbb Z})=\{0\}.$ The
same argument above also shows that $[\tau]|_{K_1(A, {\mathbb
Z}/k{\mathbb Z})}=0,$ $k=2,3,....$

\end{proof}

\begin{Remark}
{\rm It should be noted that, since $Aff(T)/\rho_B(K_0(B))$ is
divisible, $K_0(M(B)/B)/kK_0(M(B)/B)=K_0(B)/kK_0(M(B)/B).$
Therefore one sees that, for any nonzero \hm\, $\gamma: K_0(A)
\to Aff(T)/\rho_B(K_0(B)),$   there is $\alpha\in
{\rm Hom}_{\Lambda}(\underline{K}(A), \underline{K}(M(B)/B))$ such that
$\alpha|_{K_0(A)}=\gamma$
but $[\alpha]|_{K_0(A, {\mathbb Z}/k{\mathbb Z})}=0,$
$k=1,2,....$ Furthermore, if $K_1(B)$ is also divisible (or
$K_1(B)=\{0\}$), one computes that $\tau_{*1}=0$ and ${\rm
im}\tau_{*0}\subset Aff(T)/\rho_B(K_0(B))$ imply that
$$
[\tau]|_{K_0(A,{\mathbb Z}/k{\mathbb Z})}=0,\,\,\,k=2,3,...,
$$
by using the six-term exact sequence in \ref{IIDuct}.
One should note that $[\tau]|_{K_1(A, {\mathbb Z}/k{\mathbb
Z})}=0$ for $k=2,3,...,$ implies that $\tau_{*0}({\rm
tor}(K_0(A)))=0.$
On the other hand, if $\tau_{*0}({\rm tor}(K_0(A))=0$ and ${\rm
ker}\rho_B$ is divisible (or ${\rm ker}\rho_B(K_0(B))=\{0\}$),
then
$$
[\tau]|_{K_1(A, {\mathbb Z}/k{\mathbb Z})}=0,\,\,\,k=2,3,....
$$

 }
\end{Remark}

\begin{Prop}\label{IVCqa}
Let $A$ be the closure of $\cup_{n=1}^{\infty} A_n,$ where each
$A_n$ is a separable amenable \CA\, in ${\cal N}$ and let $j_n:
A_n\to A$ be the embedding. Suppose that $\tau: A\to M(B)/B$ is an
essential extension such that $\tau\circ j_n$ is  a quasidiagonal
extension for each $n.$  Then $\tau$ is also a quasidiagonal
extension.
\end{Prop}

\begin{proof}
The proof is almost the exactly the same as that of
\ref{IVLqalim}.
\end{proof}

\begin{Def}\label{IVDEm}
{\rm Denote by ${\cal C}_{afem}$ the class of separable \CA s $A$
satisfying the following: there is an embedding $j: A\to C$ such
that $j_{*0}:K_0(A)/{\rm tor}(K_0(A))\to K_0(C)$ is injective,
where $C$ is a unital AF-algebra.

Clearly every  AF-algebra is in ${\cal C}_{afem}.$ It is easy to
see that \CA s of the form $C(X)\otimes M_n$ are in ${\cal
C}_{afem},$ where $X$ is a finite CW complex. But much more is
true.

Recall that a \CA\, $A$ is called {\it residually finite
dimensional} if there is a separating family $\Pi$ of finite
dimensional irreducible representations of $A,$ i.e.,  for any
$a\in A,$ there is $\phi\in \Pi$ such that $\phi(a)\not=0.$
}\end{Def}

The following is a modification of a Dadarlat's construction.

\begin{thm}\label{IVemb}
Let $A$ be a separable amenable residually finite dimensional
\CA\, in ${\cal N}.$ Then there exists a separable unital  simple
AF-algebra $C$ and an embedding $j: A\to C$ such that $j_{*0}$
induces an injective map from $K_0(A)/{\rm tor}(K_0(A))$ into
$K_0(C).$ In particular, $A\in {\cal N}\cap {\cal C}_{afem}.$
\end{thm}

\begin{proof}
Fix  a separating sequence
of finite dimensional irreducible representations
$\{t_n\}.$ For convenience, we assume that each $t_n$
repeats infinitely many times in the sequence.
Suppose that $t_n(A)$ has rank $k(n).$
For each $n,$ define $\psi_n: A\to M_{k(n)}$ by the composition:
$A{\stackrel{t_n}{\to}} M_{k(n)}{\stackrel{\id\otimes 1_A}{\to}}
M_{k(n)}(A).$
We define a \hm\, $h_1: A\to M_{I(2)}(A),$
where $I(2)=1+k(1),$ by
$$
h_1(a)={\rm diag}(a, \psi_1(a))\,\,\,{\rm for}\,\,\, a\in A.
$$
Suppose that $h_m: M_{I(m)}(A)\to M_{I(m+1)}(A)$ is defined.
Defined  $h_{m+1}: M_{I(m+1)}(A)\to M_{I(m+2)}(A)$ by
$$
h_{m+1}(a)={\rm diag}(a, {\bar \psi}_1(a), {\bar \psi}_2(a),...,
{\bar \psi}_{m+1}(a))\,\,\,\,{\rm for}\,\,\, a\in M_{I(m+1)}(A),
$$
where $I(m+2)=I(m+1)(1+\sum_{i=1}^{m+1}k(i))$ and ${\bar
\psi}_i=\psi_i\otimes \id_{I(m+1)},$ $i=1,2,...,m+1.$ Set
$B=\lim_{m\to\infty} (M_{I(m)}(A), h_m).$ It is shown (see 3.7.8
and 3.7.9 of \cite{Lnb}) that $B$ is a unital separable amenable
simple \CA\, with $TR(B)=0.$ Since each $M_{I(m)}(A)$ satisfies
the UCT, so does $B.$ It follows from \cite{Lncl} that $B$ is
isomorphic to a unital simple AH-algebra with real rank zero and
with no dimension growth. Let $C$ be a unital simple AF-algebra
with $K_0(C)=K_0(B)/{\rm tor}(K_0(B)).$ It follows from \cite{EG}
that there exists a monomorphism $\phi$ from $B$ into $C$ such
that $\phi_{*0}$ is the quotient map from $K_0(B)$ onto
$K_0(B)/{\rm tor}(K_0(B))=K_0(C).$

Thus it remains to show that  $h_{1,\infty}: A\to B$ induces an
injective map $(h_{1,\infty})_{*0}$ on $K_0(A).$

It suffices to show that $(h_m)_{*0}$ is injective. Suppose that
$p$ and $q$ are two projections in $M_{I(m)}(A)$ such that
$h_m(p)$ and $h_m(q)$ are equivalent. Then ${\bar \psi}_j\circ
h_m(p)$ and ${\bar \psi}_j\circ h_m(q)$ are equivalent (in a
matrix algebra) for each $j.$ Therefore
$$
{\rm diag}({\bar \psi}_1(p),...,{\bar \psi}_{m+1}(p),{\bar
\psi}_1(p),...,{\bar \psi}_{m+1}(p))\andeqn {\rm diag}({\bar
\psi}_1(q),,...,{\bar \psi}_{m+1}(q),{\bar \psi}_1(q),,...,{\bar
\psi}_{m+1}(q))
$$
are equivalent. Let $t$ be any finite dimensional representation
of $M_{I(m)}(A).$ Then $(t\oplus t)\circ h_m(p)$ and $(t\oplus
t)\circ h_m(q)$ are equivalent (in a matrix algebra). From the
above, it follows that $t(p)\oplus t(p)$ and $t(q)\oplus t(q)$ are
equivalent in a matrix algebra. Thus $t(p)$ and $t(q)$ are
equivalent in the matrix algebra. This in turn implies that
$$
{\rm diag}(0, {\bar \psi}_1(p),{\bar \psi}_2(p),...,{\bar \psi}_m(p))
\andeqn
{\rm diag}(0, {\bar \psi}_1(q), {\bar \psi}_2(q),...,{\bar \psi}_m(q))
$$
are equivalent. Consequently
$$
[p]=[q] \,\,\,{\rm in}\,\,\, K_0(M_{I(m)}(A)).
$$
This implies that $(h_m)_{*0}$ is injective for each $m.$

\end{proof}

\begin{thm}\label{II1Tqd}
Let $A$ be a separable amenable \CA\, in ${\cal N}\cap {\cal
C}_{afem}.$ Suppose that $\tau$ is an essential extension.

Then $\tau$ is quasidiagonal if and only if $\tau_{*1}=0,$
 ${\rm im} \tau_{*0}\subset
Aff(T(B))/\rho(K_0(B))$ and $[\tau]|_{K_i({\mathbb Z}/k{\mathbb
Z})}=0,$ $i=0,1$ and $k=2,3,....$

\end{thm}

\begin{proof}
The ``if only" part follows from \ref{IVLqanes}.  For the ``if"
part, we first assume that $A$ is an AF-algebra. We may write
$A=\overline{\cup_{n=1}^{\infty} A_n},$ where $A_n\subset A_{n+1}$
and each $A_n$ is a finite dimensional \CA. Let $j_n: A_n\to A$ be
the embedding. It follows from \ref{IVCqa} that it suffices to
show that $\tau\circ j_n$ are quasidiagonal. Note that $(\tau\circ
j_n)_{*0} \subset Aff(T)/\rho_B(K_0(B)).$ It follows from
\ref{II1Llift} that each $\tau\circ j_n$ is in fact trivial and
therefore quasidiagonal (since $A_n$ is finite dimensional).

For the general case, let $C$ be a unital  AF-algebra and $j: A\to
C$ be an embedding such that $j_{*0}$ induces an injective \hm\,
from $K_0(A)/{\rm tor}(K_0(A))$ into $K_0(C).$ Let $\tau$ be as in
the theorem. Since $Aff(T)/\rho_B(K_0(B))$ is divisible, there
exists a \hm\, $\alpha: K_0(C)\to Aff(T)/\rho_B(K_0(B))$ such that
$\alpha\circ j_{*0}=\tau_{*0}.$ It follows from \ref{IITM} that there
is an essential extension $t: C\to M(B)/B$ such that
$t_{*0}=\alpha.$ From what we have shown, $t$ is quasidiagonal.
Let $\tau_0=t\circ j.$

Since $C$ is an AF-algebra, $K_0(C,{\mathbb Z}/k{\mathbb
Z})=K_0(C)/kK_0(C),$ $k=2,3,....$ On the other hand, ${\rm
im}t_{*0}\subset  Aff(T)/\rho_B(K_0(B))$ and
$Aff(T)/\rho_B(K_0(B))$  is divisible, so one computes that
$[t]|_{K_0(C, {\mathbb Z}/k{\mathbb Z})}=0$ for all $k,$
by using the six-term exact sequence in \ref{IIDuct}. Thus
$[\tau_0]|_{K_0(A,{\mathbb Z}/k{\mathbb Z})}=0$ for all $k.$ We
also have $(\tau_0)_{*1}=0.$

Since $C$ is an AF-algebra, $K_1(C,{\mathbb Z}/k{\mathbb
Z})=\{0\}$ for $k=2,3,....$ Since $\tau_0$ factors through $C,$ we
conclude that $[\tau_0]|{K_1(A,{\mathbb Z}/k{\mathbb Z})}=0.$
Furthermore $(\tau_0)_{*0}=\tau_{*0}.$ We then conclude that
$$
[\tau]=[\tau_0] \,\,\,{\rm in}\,\,\, KL(A, M(B)/B).
$$
Therefore, by \ref{IabTTad}, $\tau$ and $\tau_0$ are
approximately unitarily equivalent. We have shown that $t$ is a
quasidiagonal extension. So is $\tau_0,$  by \ref{IVCqa}. It
follows that $\tau$ is a quasidiagonal extension.
\end{proof}

For the last theorem in this section, one should note that
every strong NF-algebra is an inductive limit of amenable
residually finite dimensional \CA s (see 6.16 in \cite{BK}).

\begin{thm}\label{IVTqaT}
Let $A$ be the closure of $\cup_{n=1}^{\infty} A_n,$ where each
$A_n$ is a separable amenable residually finite dimensional \CA\,
in ${\cal N}.$  Let $\tau: A\to M(B)/B$ be an essential extension.
Then $\tau$ is quasidiagonal if and only if $\tau_{*1}=0,$ ${\rm
im}\tau_{*0}\subset Aff(T)/\rho_B(K_0(B))$ and $[\tau]|_{K_i(A,
{\mathbb Z}/k{\mathbb Z})}=0,$ $i=0,1$ and $k=2,3,...$
\end{thm}

\begin{proof}
It follows from \ref{IVLqanes} that we only need to prove the ``if" part
of the theorem.

 Fix an integer $n\ge 1.$ Let $\phi_n: A_n\to A$ be
the embedding. Put $\tau_n=\tau\circ \phi_n.$  So $\tau_n$ is an
essential extension. Then $(\tau_n)_{*1}=0,$ ${\rm
im}(\tau_n)_{*0}\subset Aff(T)/\rho_B(K_0(B))$  and
$[\tau_n]|_{K_i(A_n,{\mathbb Z}/k{\mathbb Z})}=0$ for $i=0,1$ and
for $k=2,3,....$
 It follows from
\ref{II1Tqd} that $\tau_n$ is quasidiagonal. Therefore the theorem
follows from \ref{IVCqa}.
\end{proof}

\section{Approximately trivial extensions}

Let
$$
0\to {\cal K}\to E\to A\to 0
$$
be an essential extension for a amenable quasidiagonal \CA\, $A.$
It is shown that the extensions is quasidiagonal if and only if it
is approximately trivial (see \cite{Br2} and \cite{Sch2}).

In this section, we will show that there are quasidiagonal extensions
that are not approximately trivial. The obstruction of a quasidiagonal
extension to be approximately trivial can be computed.
We will also discuss when an essential extension is approximately
trivial.

{\it Throughout this section $B$ is  always a non-unital but
$\sigma$-unital simple \CA\, with real rank zero, stable rank one,
weakly unperforated $K_0(B)$
and a continuous scale.}

If $A$ is a unital separable quasidiagonal \CA\, then
$A$ admits at least one tracial state. Let $T_A$ denote the tracial
state space. Let $\rho_A: K_0(A)\to Aff(T_A)$ be defined
by $\rho_A([p])(t)=t(p)$ for projection $p\in M_k(A),$ $k=1,2,...$
This map $\rho_A$ is a positive \hm. In  general,
${\rm ker}\rho_A$ is not zero.

\begin{thm}\label{IITapptr}
Let $A$ be a unital separable amenable \CA\, in ${\cal N}.$
Let $\tau:
A\to M(B)/B$ be an essential extension of $A$ by $B.$ If $\tau$ is
approximately trivial, then $\tau_{*1}=0,$ ${\rm
im}\tau_{*0}\subset Aff(T)/\rho_B(K_0(B)),$ $\tau_{*0}|_{{\rm
ker}\rho_A}=0$ and $[\tau]|_{K_i(A,{\mathbb Z}/k{\mathbb Z})}=0$
for $i=0,1$ and for $k=2,3,....$
\end{thm}

\begin{proof}
Suppose that there are monomorphisms $t_n: A\to M(B)$ such that
$$
\lim_{n\to\infty}\pi\circ t_n(a)=\tau(a)\,\,\,{\rm for\,\,\,all}\,\,\, a\in A,
$$
where $\pi: M(B)\to M(B)/B$ is the quotient map.
We obtain a positive \hm\, from $K_0(A)$ into $Aff(T).$
This implies that
$$
(t_n)_{*0}({\rm ker}\rho_A))=0.
$$
Since $K_0(B)=Aff(T)$ is a torsion free divisible group and
$K_1(M(B))=0,$ we compute that, by using the six-term exact
sequence in \ref{IIDuct},
$$
K_i(M(B),{\mathbb Z}/k{\mathbb Z})=\{0\},\,\,\, i=0,1, k=2,3,...
$$
Thus
$$
(t_n)_{*1}=0\andeqn [t_n]|_{K_i(A, {\mathbb Z}/k{\mathbb
Z})}=0,\,\,\, i=0,1,\,\,\,k=2,3,...,\,\,\,n=1,2,....
$$
For any finite subset ${\cal P}\subset {\underline{K}}(A),$ there
is an integer $n_0\ge 1$ such that
$$
[\tau_n]|_{\cal P}=[\tau]|_{\cal P}\,\,\,{\rm for\,\,\,all}\,\,\,
n\ge n_0.
$$
 Thus
 $$(\pi\circ t_n)_{*0}({\rm ker}\rho_A)=0, \,(\pi\circ
 t_n)_{*1}=0\andeqn
[\pi\circ t_n]|_{K_i(A, {\mathbb Z}/k{\mathbb
Z})}=0,\,\,\,i=0,1,\, k=2,3,...,\,n=1,2,....
$$
 Therefore
 $$
 \tau_{*0}({\rm
ker}\rho_A)=0,\,\tau_{*1}=0\andeqn [\tau]|_{K_i(A, {\mathbb
Z}/k{\mathbb Z})}=0,\,\,\,i=0,1,\, k=2,3,...,\,n=1,2,....
$$

We also have ${\rm im}(\pi\circ t_n{*0})\subset
Aff(T)/\rho_B(K_0(B)).$ It follows that
$$
{\rm im}\tau_{*0}\subset Aff(T)/\rho_B(K_0(B)).
$$
\end{proof}

\begin{Cor}\label{IICqdnotat}
Let $A$ be a unital separable AF-algebra such that
${\rm ker}\rho_A\not=0.$
Then there are quasidiagonal extensions of $A$ by $B$
that are not approximately trivial.
\end{Cor}

\begin{proof}
Let $g_0\in {\rm ker}\rho_A$ be a nonzero element.
Take any nonzero $x\in Aff(T)/\rho_B(K_0(B))$ and define a
 nonzero \hm\,
$\alpha_0: {\mathbb Z}g_0\to {\mathbb Z}x\in
Aff(T)/\rho_B(K_0(B))$ by $\alpha_0(mg_0)=mx$
for $m\in {\mathbb Z}.$ Since $ Aff(T)/\rho_B(K_0(B))$ is
divisible, we obtain a \hm\, $\alpha: K_0(A)\to
Aff(T)/\rho_B(K_0(B))$ such that $\alpha(g_0)=\alpha_0(g_0).$ It
follows from \ref{IIILemb} that there is a monomorphism $\tau:
A\to M(B)/B$ such that $\tau_{*0}=\alpha.$ It follows from
\ref{II1Tqd} that $\tau$ is quasidiagonal. But by \ref{IITapptr}
$\tau$ is not approximately trivial.
\end{proof}

\begin{Remark}
{\rm From \ref{IICqdnotat}, one sees that it is typical rather than unusual
that quasidiagonal extensions are different from approximately trivial
extensions. The assumption that $A$ is AF is certainly not necessary.
}
\end{Remark}

\begin{Lem}\label{IVLgroup}
Let $G$ be an unperforated ordered group with the Riesz interpolation
property. Suppose that $G_0\subset G$ is a countable ordered
subgroup. Then there exists a countable subgroup $G_1\supset G_0$ which
satisfies the Riesz interpolation property and is unperforated.
If $G$ is simple, we may further assume that $G_1$ is also simple.
\end{Lem}

\begin{proof}
Since $G_0$ is countable, there exists a countable ordered subgroup
$F_1$ of $G$ such that $G_0\subset F_1$ and
if $g_1, g_2, g_3, g_4\in G_0$ with $g_1, g_2\le g_3, g_4$
then there is $g\in F_1$ such that
$$
g_1, g_2\le h \le g_3, g_4.
$$
If a countable ordered subgroup $F_n$ has been constructed,
we have a countable ordered subgroup $F_{n+1}$ such that
if $x_1,x_2\le y_1,y_2$ are in $F_n$ there exists
$g\in F_{n+1}$ such that
$$
x_1, x_2\le g\le y_1, y_2.
$$
Set $G_1=\cup_{n=1}^{\infty} F_n.$
Then $G_1$ is a countable ordered subgroup of $G$ containing
$G_0.$ From the construction, it is also clear
that $G_1$ has the Riesz interpolation property.

Now we further assume that $G$ is simple. Let $g\in (G_1)_+$ be a nonzero
positive element and $f\in G_1.$ Since $G$ is simple there is an
integer $n\ge 1$ such that $ng\ge f.$ This implies that $G_1$ is
also simple. Since $G$ is unperforated, it follows that $G_1$ is
also unperforated.
\end{proof}

\begin{Lem}\label{IVLabemb}
Let $A$ be a unital separable commutative \CA. Then there exists a
monomorphism $h: A\to M(B)$ such that $\pi\circ h$ is an essential
extension and ${\rm im} h_{*0}\subset\rho_B(K_0(B)).$
\end{Lem}

\begin{proof}
Let $G\subset \rho_B(K_0(B))$ be a countable simple ordered group
with the Riesz interpolation property. We may also assume that
$G\not\cong {\mathbb Z}.$ It follows \cite{EHS} that there is a
unital non-elementary simple AF-algebra $C$ such that $K_0(C)$ is
order isomorphic to $G.$ Thus we obtain an order isomorphism
$\alpha: K_0(C)\to G\subset \rho_B(K_0(B))\subset Aff(T).$ It
follows from p. 67 in \cite{AS} that there is monomorphism $j: A\to C.$
It follows from \ref{IIILemb} that there is a \hm\, $\phi: C\to
M(B)$ such that $\phi(A)\cap B=\{0\}$ and $\phi_{*0}=\alpha.$
Define $h=\phi\circ j.$
\end{proof}

\begin{thm}\label{IVTabeappr}
Let $A$ be a separable unital commutative \CA. Then $\tau: A\to
M(B)/B$ is approximately trivial if and only if $\tau_{*1}=0,$
$\tau_{*0}({\rm ker}\rho_A)=0,$ ${\rm im}\tau_{*0}\subset
Aff(T)/\rho_B(K_0(B))$ and $[\tau]|_{K_i(A,{\mathbb Z}/k{\mathbb
Z})}=0,$ $i=0,1,$ $k=2,3,....$
\end{thm}

\begin{proof}
We may write $A=C(X),$ where $X$ is a compact metric space. There
are finite CW complexes $X_n$ such that $X=\lim_{\leftarrow n}
X_n.$ We also write $C(X)=\lim_{n\to\infty} C(X_n)$ and $\psi_n:
C(X_n)\to C(X)$ is the induced \hm. We will show that if
$\tau_{*1}=0,$ $\tau_{*0}({\rm ker}\rho_A)=0,$ ${\rm
im}\tau_{*0}\subset Aff(T)/\rho_B(K_0(B))$ and
$[\tau]|_{K_i(A,{\mathbb Z}/k{\mathbb Z})}=0,$ $i=0,1$ and
$k=2,3,...,$ then $\tau$ is approximately trivial. Let
$e_1^{(n)},e_2^{(n)},...,e_{k(n)}^{(n)}$ be projections
corresponding to each summand of $C(X_n)$ which corresponds to
each connected component of $X_n.$ So $K_0(C(X_n))/{\rm
ker}\rho_C$ is generated by
$\{e_1^{(n)},e_2^{(n)},...,e_{k(n)}^{(n)}\}.$ Denote by $G_n$ the
subgroup of $K_0(C(X))$ generated by $[\psi_n(e_1^{(n)})],
[\psi_n(e_2^{(n)})], ..., [\psi_n(e_{k(n)}^{(n)})].$ Let $f_n:
X\to X_n$ denote the continuous map induced by the inverse inductive
limit $X=\lim_{\leftarrow n} X_n.$ Let
$\xi_1^{(n)},\xi_2^{(n)},...,\xi_{k(n)}^{(n)}$ be points in $X_n$
which lie in different components. Let
$x_1^{(n)},x_2^{(n)},...,x_{k(n)}^{(n)}$ be points in $X$ such
that $f_n(x_i^{(n)})=\xi_i^{(n)}.$ Since $\rho_B(K_0(B))$ is dense
in $Aff(T),$ we can find mutually orthogonal projections
$p_1^{(n)},p_2^{(n)},..., p_{k(n)}^{(n)}$ in $M(B)$ such that
$\sum_{i=1}^{k(n)}p_i^{(n)}\le  1_{M(B)},$
$1_{M(B)}-\sum_{i=1}^{k(n)}p_i^{(n)}\not\in B$ and
$[\pi(p_i^{(n)})]=[\tau\circ \psi_n(p_i^{(n)})],$
$i=1,2,...,k(n).$ Define a \hm\, $\phi_n: C(X)\to M(B)$ by
$$
\phi_n(f)=\sum_{i=1}^{k(n)}f(x_i^{(n)})p_i^{(n)}\,\,\,{\rm for}\,\,\,
f\in C(X).
$$
Let $P_n\le 1_{M(B)}-\sum_{i=1}^{k(n)}p_i^{(n)}$ be a projection
in $M(B)\setminus B$ such that $[P_n]\in \rho_B(K_0(B)).$ It
follows from \ref{IVLabemb} that there is a monomorphism
$h_n^{(0)}: C(X)\to P_nM(B)P_n$ such that ${\rm
im}(h_n^{(0)})_{*0} \subset \rho_B(K_0(B))$ and $\pi\circ
h_n^{(0)}$ is injective. Now define $h_n=\phi_n+h_n^{(0)}.$ Then
$h_n: C(X)\to M(B)$ give an essential trivial extension of $C(X)$
by $B.$ Note that $(\pi\circ h_n)_{*1}=0,$ $[\pi\circ
h_n]|_{K_i({\mathbb Z}/k{\mathbb Z})}=0,$ $i=0,1,$  $k=2,3,...,$
and
$$
(\pi\circ h_n)_{*0}|_{G_n}=\tau_{*0}|_{G_n}, n=1,2,...,
$$
and ${\rm im}(\pi\circ h_n)_{*0}) \subset Aff(T)/\rho_B(K_0(B)).$
Since $\cup_{n=1}^{\infty} G_n=K_0(C(X)),$
 one checks that, for any
finite subset ${\cal P}\subset {\underline{K}}(C(X)),$ there
exists an integer $n,$ such that
$$
[\pi\circ h_n]|_{\cal P}=[\tau]_{\cal P}.
$$
It follows from \ref{IabT1} that there exists a sequence of
unitaries $u_n\in M(B)/B$ such that
$$
\lim_{n\to\infty} {\rm ad}\,u_n\circ \pi\circ
h_n(a)=\tau(a)\,\,\,a\in C(X).
$$
This implies that $\tau$ is approximately trivial.
\end{proof}

\begin{Def}\label{IIDpl}
{\rm Let
$
Pl(K_0(A), Aff(T)/\rho_B(K_0(B)))
$
be the set of  those elements $\alpha$ in $ {\rm Hom}(K_0(A), Aff(T)/\rho_B(K_0(B))$ such that
there exists a positive \hm\, $\beta: K_0(A)\to Aff(T)$ such that
$\Phi\circ \beta=\alpha,$
where  $\Phi: Aff(T)\to Aff(T)/\rho_B(K_0(B))$ is the quotient map.

Denote by $ Apl(K_0(A), Aff(T)/\rho_B(K_0(B))) $  the set of those elements
$\alpha $ in $ {\rm Hom}(K_0(A), Aff(T)/\rho_B(K_0(B))$ satisfying the
following: there exists an increasing  sequence of finitely generated subgroups
$\{G_n\}\subset K_0(A)$ such that $\cup_{n=1}^{\infty} G_n=K_0(A)$
and a sequence of \hm\, $\alpha_n\in Pl(K_0(A),
Aff(T)/\rho_B(K_0(B)))$ such that
$$
(\alpha_n)|_{G_n}=\alpha|_{G_n},
$$
$n=1,2,....$

One should note that if $\alpha\in Apl(K_0(A). Aff(T)/\rho_B(K_0(B))$
then $\alpha|_{{\rm ker}\rho_A}=0.$
In fact, if $x\in {\rm ker}\rho_A,$ then $x\in G_n$ for some $n.$
But $\alpha_n(x)=0$ for all $n.$ Therefore $\alpha(x)=0.$
}
\end{Def}

\begin{Prop}\label{IITnappext}
Let $A$ be a separable amenable \CA\, and let
$\tau: A\to M(B)/B$ be an essential approximately
trivial extension.
Then
$$
\tau_{*1}=0, [\tau]|_{K_i(A,{\mathbb Z}/k{\mathbb Z})}=0\andeqn
$$
$$\tau_{*0}\in Apl(K_0(A), Aff(T)/\rho_B(K_0(B))).$$
\end{Prop}

\begin{proof}
Suppose that $t_n: A\to M(B)/B$ be a sequence of trivial
extensions such that
$$
\lim_{n\to\infty}t_n(a)=\tau(a)\,\,\,\,{\rm for\,\,\,all}\,\,\, a\in A.
$$
There is a sequence of monomorphism $h_n: A\to M(B)$ such that
$\pi\circ h_n=t_n,$ where $\pi: M(B)\to M(B)/B$ is the quotient
map, $n=1,2,....$ Since $K_1(M(B))=0$ and $K_0(M(B))$ has no
torsion, we conclude (by using the six-term exact sequence in
\ref{IIDuct}) that
$$
[t_n]|_{K_1(A)}=0\andeqn [t_n]|_{K_i(A, {\mathbb Z}/k{\mathbb Z})}=0
\,\,\,i=0,1,\, k=1,2,....$$
It follows that
$$
[\tau]|_{K_1(A)}=0 \andeqn [t]|_{K_i(A, {\mathbb Z}/k{\mathbb Z})}=0,
\,\,\,i=0,1,\, k=1,2,....$$
Let $\{G_n\}$ be a sequence of finitely generated groups
of $K_0(A)$ such that $\cup_{n=1}^{\infty} G_n=K_0(A).$
For each $n,$ there is $m(n)$ such that
$$
(t_m)|_{G_n}=\tau|_{G_n}
$$
for all $m\ge m(n),$ since
$\lim_{n\to\infty}t_n(a)=\tau(a)$ for all $a\in A.$
However, $(t_m)_{*0}\in Pl(K_0(A), Aff(T)/\rho_B(K_0(B))).$
Therefore $\tau_{*0}\in Apl(K_0(A), Aff(T)/\rho_B(K_0(B))).$
\end{proof}

\begin{thm}\label{IITappext}
Let $A$ be a  separable amenable \CA\, in ${\cal N}$
which can be embedded into a unital AF-algebra $C$ such that
$K_0(A)/{\rm ker}\rho_A=K_0(C)/{\rm ker}\rho_C$ (as ordered
groups).
Let $\tau: A\to M(B)/B$ be an essential extension.
Then $\tau$ is approximately trivial if and only if
$\tau_{*1}=0,$ $[\tau]|_{K_i(A, {\mathbb Z}/k{\mathbb Z})}=0,$ $i=0,1,$
$k=1,2,...,$ and $\tau_{*0}\in Apl(K_0(A), Aff(T)/\rho_B(K_0(B))).$
\end{thm}

\begin{proof}
The ``onlyi if part" follows from \ref{IITnappext}.
Suppose that $\alpha_n\in Pl(K_0(A), Aff(T)/\rho_B(K_0(B)))$
is such that
$$
(\alpha_n)|_{G_n}=\tau_{*0}|_{G_n},
$$
where $\cup_{n=1}^{\infty}G_n=K_0(A)$ and $G_n$ is finitely generated.
Suppose
also $h_n: K_0(A)\to K_0(M(B))=Aff(T)$ is a positive \hm\, such that
$\Phi\circ h_n=\alpha_n.$ Note that $h_n|_{{\rm ker}\rho_A}=0.$
Let ${\tilde h}_n: K_0(A)/{\rm ker}\rho_A\to K_0(M(B))$ denote
the induced positive \hm.
Let $C$ be the unital simple AF-algebra
such that there exists an embedding
$j: A\to C$ such that $j_{*0}$ induces an order isomorphism
from $K_0(A)/{\rm ker}\rho_A$ onto $K_0(C)/{\rm ker}\rho_C$
with $K_0(C)=K_0(A)/{\rm ker}\rho_A$
and $[1_C]=[1_A].$
There is a \hm\, $\psi_n: C\to M(B)$ such that
$(\psi_n)_{*0}={\tilde h_n}.$
Put $\phi_n=\psi_n\circ j.$
Let $t_n=\pi\circ \phi_n.$
It is clear that $[t_n]_{K_1(A)}=0$
and $(t_n)_{*0}=\alpha_n.$
We note that since $K_0(M(B))=Aff(T)$ is divisible,
$K_0(M(B))/kK_0(M(B))=0.$  We also have $K_1(M(B))=0.$
This implies that $K_0(M(B),{\mathbb Z}/k{\mathbb Z})=0$ for all
$k.$  Therefore
$$
[t_n]|_{K_0(A, {\mathbb Z}/k{\mathbb Z})}=0\,\,\,\, k=1,2,....
$$
Since
$(\phi_n)_{*1}=0$ and $K_0(M(B))$ has no torsion,
we compute that
$$
[t_n]|_{K_1(A)}=0\andeqn [t_n]|_{K_1(A, {\mathbb Z}/k{\mathbb Z})}=0
$$
for all $k.$ For  $\tau,$ we note that $\tau_{*0}\in APl(K_0(A),
Aff(T)/\rho_B(K_0(B)))$ implies that $\tau_{*0}|_{{\rm
ker}\rho_A}=0.$ In particular, $\tau_{*0}|_{{\rm tor}(K_0(A))}=0.$
Since $\tau_{*1}=0,$ it follows that
$$
[\tau]|_{K_i(A, {\mathbb Z}/k{\mathbb Z})}=0\,\,\,i=0,1\,\,\,k=1,2,....
$$
Therefore, for any finite subset ${\cal P}\subset
{\underline{K}(A)},$ there exists $n$ such that
$$
[\tau]|_{{\cal P}}=[t_m]|_{{\cal P}}
$$
for all $m\ge n.$ It follows from \ref{IabT1} that there are
unitaries $u_n\in M(B)/B$ such that
$$
\lim_{n\to\infty} {\rm ad}\, u_n\circ  t_n(a)=\tau(a)
$$
for all $a\in A.$

\end{proof}

\begin{Ex}\label{IIE5}
{\rm There are many examples where a \hm\, $\alpha: K_0(A)\to
Aff(T)/\rho_B(K_0(B))$ can not be lifted to a positive \hm\,
$\beta: K_0(A)\to Aff(T).$ In \ref{IIe2}, an example is given where
even if $\alpha=0,$ one can not get a nonzero positive \hm\,
$\beta: K_0(A)\to Aff(T)$ such that $\Phi\circ \beta=\alpha.$
To show other complications, let
$\rho_B(K_0(B))={\mathbb Z}[1/2]$ and let $K_0(A)={\mathbb
Z}\oplus {\mathbb Z}\sqrt{2}.$ Define a \hm\, $\alpha: K_0(A)\to
Aff(T)/\rho_B(K_0(B))= {\mathbb R}/{\mathbb Z}[1/2]$ so that
$\alpha(1)=1$ and $\alpha(\sqrt{2}) =\pi.$ If
there is  a nonzero \hm\, $\beta: K_0(A)\to {\mathbb R}$ such that
$\Phi\circ \beta(1)=\alpha(1)$ and $\Phi\circ
\beta(\sqrt{2})=\alpha(\sqrt{2}),$ then $\beta(1)=x$ and
$\beta(\sqrt{2})=y$ with $x\in {\mathbb Z}[1/2]$ and $y=\pi+z,$
where $z\in {\mathbb Z}[1/2].$ If, in addition, $\beta$ were
positive, then $\beta(\sqrt{2})=x\sqrt{2}.$ But that would imply
that $x\sqrt{2}=\pi+z.$ This is impossible since $x, z\in {\mathbb Z}[1/2].$
 Therefore $\alpha\not\in Apl(K_0(A),
Aff(T)/\rho_B(K_0(B))).$

This example also shows that there are very few elements in
$APl(K_0(A), Aff(T)/\rho_B(K_0(B)))$ or in \\
$Pl(K_0(A), Aff(T)/\rho_B(K_0(B))).$

Nevertheless, we have the following:
}
\end{Ex}

\begin{thm}\label{IIP1}
Let $A$ be a unital separable \CA\, with $K_0(A)={\mathbb Z}[1/p],$ where
$p$ is a prime number.
 Suppose that $\tau_{*0}:K_0(A)
\to Aff(T)/\rho_B(K_0(B))$ is a \hm. Then the following hold:

{\rm (1)} $\tau_{*0} \in Apl(K_0(A), Aff(T)/\rho_B(K_0(B))$ if
$\tau_{*0}\not=0,$

{\rm (2)} if $\rho_B(K_0(B))$ is divisible by $p,$
then every such nonzero $\tau_{*0}$ is in $Pl(K_0(A), Aff(T)/\rho_B(K_0(B)).$

{\rm (3)} if $\rho_B(K_0(B))$ is finitely generated (such as
${\mathbb Z}\oplus {\mathbb Z}_{\theta}$ for some irrational
number $\theta$), then
$$
Apl(K_0(A), Aff(T)/\rho_B(K_0(B))\not=Pl(K_0(A),Aff(T)/\rho_B(K_0(B)),
$$

{\rm (4)} $\tau_{*0}=0$ is in $Apl(K_0(A), Aff(T)/\rho_B(K_0(B)))$
if and only if there is a
sequence of nonzero positive elements $ \eta_n\in \rho_B(K_0(B))$ such that
$\eta_n\le {\widehat{1_{M(B)}}}$ and
$\eta_n$ is divisible by $p^n,$ $n=1,2,....$

\end{thm}

\begin{proof}

To prove (1), let
$\tau_{*0}: K_0(A)\to Aff(T)/\rho_B(K_0(B))$ be a nonzero \hm.
There exists $N$ such that $\xi_n=\tau_{*0}(1/p^n)\not=0$
for all $n\ge N.$
Since $\rho_B(K_0(B))$ is dense in $Aff(T),$
one can choose $r_n\gg 0$ in $Aff(T)$ such that $\Phi(r_n)=\tau_{*0}(1/p^n),$
where $\Phi: Aff(T)\to Aff(T)/\rho_B(K_0(B))$ is the quotient map.
Since $\rho_B(K_0(B))$ is dense, it is easy to find
$r_n\gg 0$ in $Aff(T)$ such that $p^nr_n\le 1_{M(B)}.$
Define $\beta_n(z)=p^nr_nz$ for $z\in K_0(A)$ ($n\ge N$). Here
we identify $z$ with the real number $z$ (so $rz\in Aff(T)$).
Write
$$
G_n=\{m/p^{n+N}: m\in {\mathbb Z}\},\,\,\,n=1,2,....
$$
Then $G_n$ is finitely generated and $\bigcup_{n=1}^{\infty} G_n=K_0(A).$
Moreover,
$$
(\Phi\circ \beta_n)|_{G_n}=\tau_{*0}|_{G_n},\,\,\, n=1,2,....
$$
This proves (1).

To prove (2), we first note  that $Aff(T)/\rho_B(K_0(B))$ has no $p$-torsion.
   Suppose $x\in Aff(T)/\rho_B(K_0(B))$ is a nonzero element
   so that $px=0.$ Let $y\in Aff(T)$ so that $\Phi(y)=x.$
   Then $py\in \rho_B(K_0(B)).$
   Since $\rho_B(K_0(B))$ is divisible by $p,$
   there is $z\in \rho_B(K_0(B))$ such that
   $p(y-z)=0.$ Since $Aff(T)$ is torsion free, $y=z,$ or
   $x=0.$

We will show that $\Phi(p^N_Nz)=\tau_{*0}(z)$ for $z\in K_0(A)={\mathbb Z}[1/p].$
Note we have shown above that for any $z\in {\mathbb Z}(1/p^N),$
$\Phi(p^Nr_Nz)=\tau_{*0}(z).$
Suppose that $x=\tau_{*0}(1/p^{N+1}).$ Then $px=\tau_{*0}(1/p^N)
=\Phi(r_N).$ Thus
$px=p\Phi\circ\beta_N(1/p^{N+1}).$ This implies
that $p(x-\xi_N)=0$ in $Aff(T)/\rho_B(K_0(B)).$
Since $Aff(T)/\rho_B(K_0(B))$ has no $p$-torsion, $x=\xi_N.$
Therefore $\Phi\circ \beta_N(z)=\tau_{*0}(z)$ for
all $z\in {\mathbb Z}(1/p^{N+1}).$ By induction,
we verify that $\Phi\circ \beta_N=\tau_{*0}.$
It is clear that $\beta_N$ is positive.

For (3),
we note that $ext_{\mathbb Z}({\mathbb Z}[1/p], {\mathbb Z})\not=\{0\}$
(see Theorem 99.1 and p. 179 and in \cite{Fu}, and  also
\cite{Rot}). Therefore, since $\rho_B(K_0(B))$ is a finite sum
of ${\mathbb Z},$ we conclude that
$ext_{\mathbb Z}({\mathbb Z}[1/p], \rho_B(K_0(B)))\not=\{0\}.$
Consider the following exact sequence:
$$
\cdots {\rm Hom}({\mathbb Z}[1/p], Aff(T))\to {\rm Hom}({\mathbb Z}[1/p], Aff(T)/\rho_B(K_0(B)))
\to ext({\mathbb Z}[1/p], \rho_B(K_0(B)))\to
$$
$$
 ext({\mathbb Z}[1/p],
Aff(T))\to ext({\mathbb Z}[1/p], Aff(T)/\rho_B(K_0(B)))\to\cdots .
$$
Since $Aff(T)$ is divisible, $ext({\mathbb Z}[1/p],Aff(T))=\{0\}.$
This implies that
the map from \\
${\rm Hom}({\mathbb Z}[1/p], Aff(T)/\rho_B(K_0(B)))$
to $ext({\mathbb Z}[1/p], \rho_B(K_0(B)))$ is surjective.
Thus we obtain $\tau_{*0}: {\mathbb Z}[1/p]\to Aff(T)/\rho_B(K_0(B))$ such
that it gives a non-splitting extension
of ${\mathbb Z}[1/p]$ by $\rho_B(K_0(B)).$
This $\tau_{*0}$ can not be lifted to a \hm\, from $K_0(A)$ into
$Aff(T).$ In particular $\tau_{*0}\not\in
Pl(K_0(A), Aff(T)/\rho_B(K_0(B)).$ It follows from (1) that $\tau_{*0}\in
Apl(K_0(A), Aff(T)/\rho_B(K_0(B)).$ Thus (3) follows.

To see  (4), suppose that
there exists a
sequence of nonzero positive elements $\eta_n\in \rho_B(K_0(B))$
such that $\eta_n\le {\widehat{1_{M(B)}}}$ and $\eta_n=p^n\zeta_n$ for
some $\zeta_n\in \rho_B(K_0(B)),$ $n=1,2,....$
Since $\rho_B(K_0(B))$ is weakly unperforated, $\zeta_n\gg 0.$
Fix $n,$
define $\beta_n: {\mathbb Z}[1/p]\to Aff(T)$ by
$\beta_n(z)=p^nz \zeta_n$ for all $z\in {\mathbb Z}[1/p]\,{\rm (}
\subset {\mathbb R} {\rm )}.$
Then $\beta_n$ is a positive \hm .
Let $G_n$ be as in the proof of (1).
Then $\beta_n(G_n)\subset \rho_B(K_0(B)),$ $n=1,2,....$
So $(\Phi\circ \beta_n)|_{G_n}=0,$ $n=1,2,....$
Therefore $\tau_{*0}=0$ is in $Apl(K_0(A), Aff(T)/\rho_B(K_0(B))).$

Conversely, if $\tau_{*0}=0$ is in  $Apl(K_0(A), Aff(T)/\rho_B(K_0(B))).$
Suppose that $F_n\subset {\mathbb Z}[1/p]$ is an increasing
sequence of finitely
generated subgroups such that $\bigcup_{n=1}^{\infty}F_n={\mathbb Z}[1/p]$
and there is a sequence of positive \hm s $\alpha_n: {\mathbb Z}[1/p]
\to Aff(T)$ such that
$\alpha_n(F_n)\subset \rho_B(K_0(B)),$ $n=1,2,....$
Replacing $\alpha_n$ by $t_n\alpha_n$ for some positive real numbers $t_n,$
we may assume that $\alpha_n(1)\le {\widehat{1_{M(B)}}},$ $n=1,2,....$
Without loss of generality, we may also assume that
$1, 1/p^n\subset F_n,$ $n=1,2,....$
Let $\eta_n=\alpha_n(1).$ Since $\alpha_n(1/p^n)\subset \rho_B(K_0(B)),$
$\eta_n$ is divisible by $p^n,$ $n=1,2,....$
\end{proof}

\begin{Remark}
Note that $\rho_B(K_0(B))$ may not have any nonzero elements
to be divided by $p.$ In these cases, the condition in (4) in the previous theorem and the next corollary never holds. In other words, when $[\tau]=0$
in $KL(A, M(B)/B),$ $\tau$ is never approximately trivial.
\end{Remark}

\begin{Cor}\label{IIICl}
Let $A$ be a unital separable amenable \CA\, in ${\cal N}.$
Suppose that there is a monomorphism $j:A \to C$ for some
 unital simple AF-algebra
$C$ with $C={\mathbb Z}[1/p].$ Suppose also
that $j_{*0}$ maps
$K_0(A)/{\rm ker}\rho_A$ injectively to ${\mathbb Z}[1/p].$
 Let  $\tau: A\to M(B)/B$ be an essential extension.

{\rm (1)} If  $[\tau]\not=0,$ $\tau_{*1}=0,$ $\tau_{*0}|_{{\rm ker}\rho_A}=0,$
 ${\rm im}\tau_{*0}\subset Aff(T)/\rho_B(K_0(B))$ and
$[\tau]|_{K_i(A,{\mathbb Z}/k{\mathbb Z})}=0,$ $i=0,1$ and $k=2,3,...$
Then $\tau$ is approximately trivial.

{\rm (2)}  If there is no positive \hm\,
from $K_0(A)$ into $\rho_B(K_0(B)),$ then
no essential extension with $[\tau]=0$ in $KL(A, M(B)/B)$ can be trivial.
Furthermore, if $A=C,$ then there exists an essential trivial extension
$\tau$ with $[\tau]=0$ in $KL(A, M(B)/B)$
if and only if there is a positive \hm\,
$\alpha: K_0(A)\to \rho_B(K_0(B)).$

{\rm (3)}  Suppose further that $K_0(A)/{\rm ker}\rho_A={\mathbb Z}[1/p].$
Then $[\tau]=0$ implies that $\tau$ is approximately
trivial if and only if there exists a
sequence of nonzero positive elements $\eta_n\in \rho_B(K_0(B))$
such that $\eta_n\le {\widehat{1_{M(B)}}}$ and $\eta_n$ is divisible
by $p^n,$ $n=1,2,....$

\end{Cor}

\begin{proof}
Suppose $j: A\to C$ is the embedding.
Since $Aff(T)/\rho_B(K_0(B))$ is divisible, there exists a \hm\,
$\alpha: K_0(C)\to Aff(T)/\rho_B(K_0(B))$ such
that $\alpha|_{j_{*0}(K_0(A))}=\tau_{*0}$ which is nonzero.
Let $\psi: C\to M(B)/B$ be an essential extension
such that $\psi_{*0}=\alpha.$
Since $Aff(T)/\rho_B(K_0(B))$ is divisible and $K_1(C)=\{0\},$
 one computes
that $[\tau]|_{K_0(C, {\mathbb Z}/k{\mathbb Z})}=0$
for all $k.$
Since $K_0(C)$ has no torsion and $K_1(C)=0,$
$[\tau]|_{K_1(C, {\mathbb Z}/k{\mathbb Z})}=0$ for all $k.$
It follows
that
$$
[\tau]=[\psi\circ j]\,\,\, {\rm in}\,\,\, KL(A, M(B)/B).
$$
So $\tau$ and $\psi\circ j$ are  approximately
unitarily equivalent. Thus it suffices to show that $\psi\circ j$ is
approximately trivial. It follows from \ref{IIP1} that
$\psi_{*0}\subset Apl(K_0(C), Aff(T)/\rho_B(K_0(B)).$
It follows  from \ref{IITappext} that $\psi$ is approximately trivial.
Hence $\psi\circ j$ is approximately trivial.

For (2),  we assume that $[\tau]=0.$
Suppose that $\tau$ is trivial and $h: A\to M(B)$ is a monomorphism
such that $\pi\circ h=\tau.$
Then $h_{*0}$ gives a positive \hm\, from $K_0(A)$ into $\rho_B(K_0(B)).$

Suppose now that $A=C.$
If there is a $\alpha: K_0(A)\to \rho_B(K_0(B))\subset Aff(T),$ then by
\ref{IIILemb} there exists a monomorphism\, $h:A\to M(B)$ such that
$h_{*0}=\alpha$ and $h(A)\cap B=\{0\}.$ Thus $\tau=\pi\circ h$ is trivial
and $[\tau]=0.$

Now consider (3). Suppose that $[\tau]=0$ in $KL(A, M(B)/B)$ and
suppose also that there is a \hm\, $h_n: A\to M(B)$ such that
$$
\lim_{n\to\infty}\pi\circ h_n(a)=\tau(a)\,\,\,{\rm for\,\,\, all}\,\,\, a\in A.
$$
Then $\tau_{*0}\in Apl(K_0(A), Aff(T)/\rho_B(K_0(B))).$
Thus the ``if only" part follows from (4) in \ref{IIP1}.
On the other hand, if those $\eta_n$ exists, by (4) in \ref{IIP1},
$\tau_{*0}\in Apl(K_0(A), Aff(T)/\rho_B(K_0(B))).$ Thus
the ``if" part follows from \ref{IITappext}.
\end{proof}

\end{document}